\documentclass[10pt]{article}

\usepackage{graphicx, amsfonts,amssymb, 
	amsmath,amsthm,color,mathtools,tikz,tikz-3dplot,hyperref,url}

\usepackage{epstopdf}
\usepackage{wasysym}
\usetikzlibrary{calc,intersections,through,backgrounds, arrows,decorations.markings,arrows.meta}
\usepackage{tkz-euclide}

\newtheorem{lemma}{Lemma}[section]
\newtheorem{proposition}[lemma]{Proposition}

\newtheorem{theorem}{Theorem}
\newtheorem{remark}[lemma]{Remark}
\newtheorem{example}[lemma]{Example}

\newtheorem{corollary}[lemma]{Corollary}

\graphicspath{{fig/}}

\begin{document}

\def\mainfile{}
\newcommand{\eps}{{\varepsilon}}
\newcommand{\C}{{\mathbb C}}
\newcommand{\Q}{{\mathbb Q}}
\newcommand{\R}{{\mathbb R}}
\newcommand{\Z}{{\mathbb Z}}
\newcommand{\RP}{{\mathbb {RP}}}
\newcommand{\CP}{{\mathbb {CP}}}
\newcommand{\Tr}{\rm Tr}
\newcommand{\g}{\gamma}
\newcommand{\G}{\Gamma}
\newcommand{\e}{\varepsilon}
\newcommand{\crel}{\stackrel{c}{\sim}}
\newcommand{\drel}{\stackrel{d}{\sim}}
\newcommand{\ed}{s}
\newcommand{\ve}{v}
\newcommand{\rec}{\mathsf{R}}

\title{A family of integrable transformations of centroaffine polygons: geometrical aspects}

\author{Maxim Arnold\footnote{
Department of Mathematics, 
University of Texas, 
800 West Campbell Road,
Richardson, TX 75080;
maxim.arnold@utdallas.edu}
\and
Dmitry Fuchs\footnote{
Department of Mathematics, 
University of California, 
Davis, CA 95616;
 fuchs@math.ucdavis.edu}
\and
Serge Tabachnikov\footnote{
Department of Mathematics,
Pennsylvania State University,
University Park, PA 16802;
tabachni@math.psu.edu}}

\date{}
\maketitle

\begin{abstract}
Two polygons, $(P_1,\ldots,P_n)$ and $(Q_1,\ldots,Q_n)$ in ${\mathbb R}^2$  are $c$-related if 
$\det(P_i, P_{i+1})=\det(Q_i, Q_{i+1})$  and $\det(P_i, Q_i)=c$ for all $i$. This relation extends to twisted polygons (polygons with monodromy), and it  descends to the moduli space of $SL(2,{\mathbb R}^2)$-equivalent polygons. This relation is an equiaffine analog of the discrete bicycle correspondence studied by a number of authors. We study the geometry of this relations, present its integrals, and show that, in an appropriate sense,  these relations, considered for different values of the constants $c$, commute. We relate this topic with the dressing chain of Veselov and Shabat. The case of small-gons is 
investigated in detail.
\end{abstract}

\tableofcontents

\section{Introduction} \label{sec:intro}

The motivation for this paper is two-fold. 

The first one is the study of the {\it discrete bicycle correspondence}  on polygons in the Euclidean plane \cite{TT}, a discretization of the {\it bicycle correspondence} on smooth curves that was studied in \cite{BLPT,Ta17}. See also \cite{Ho,PSW} for the discrete and \cite{We} for the continuous versions of this correspondence. 

\begin{figure}[ht]
\centering
\includegraphics[width=.38\textwidth]{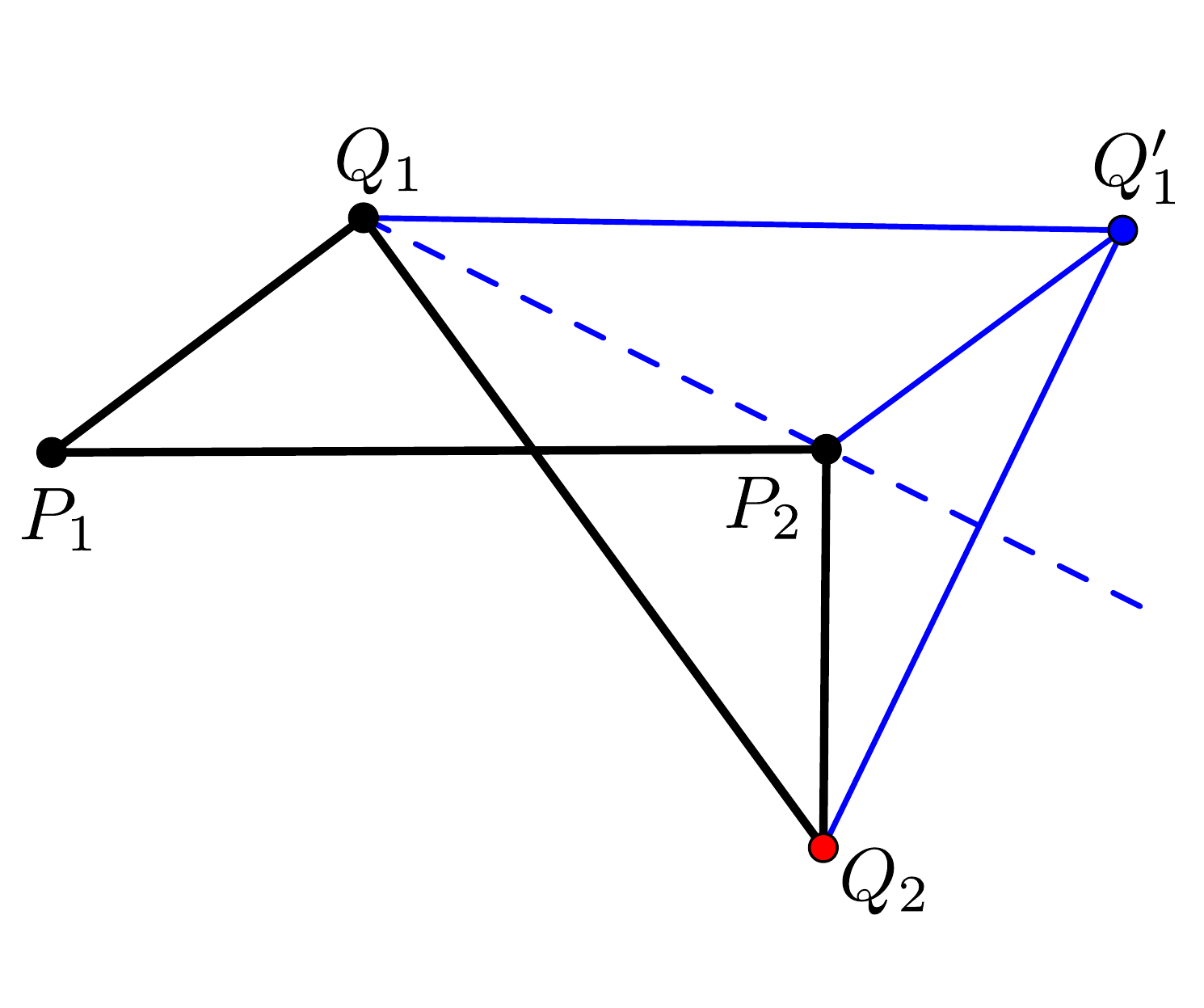}\quad
\includegraphics[width=.57\textwidth]{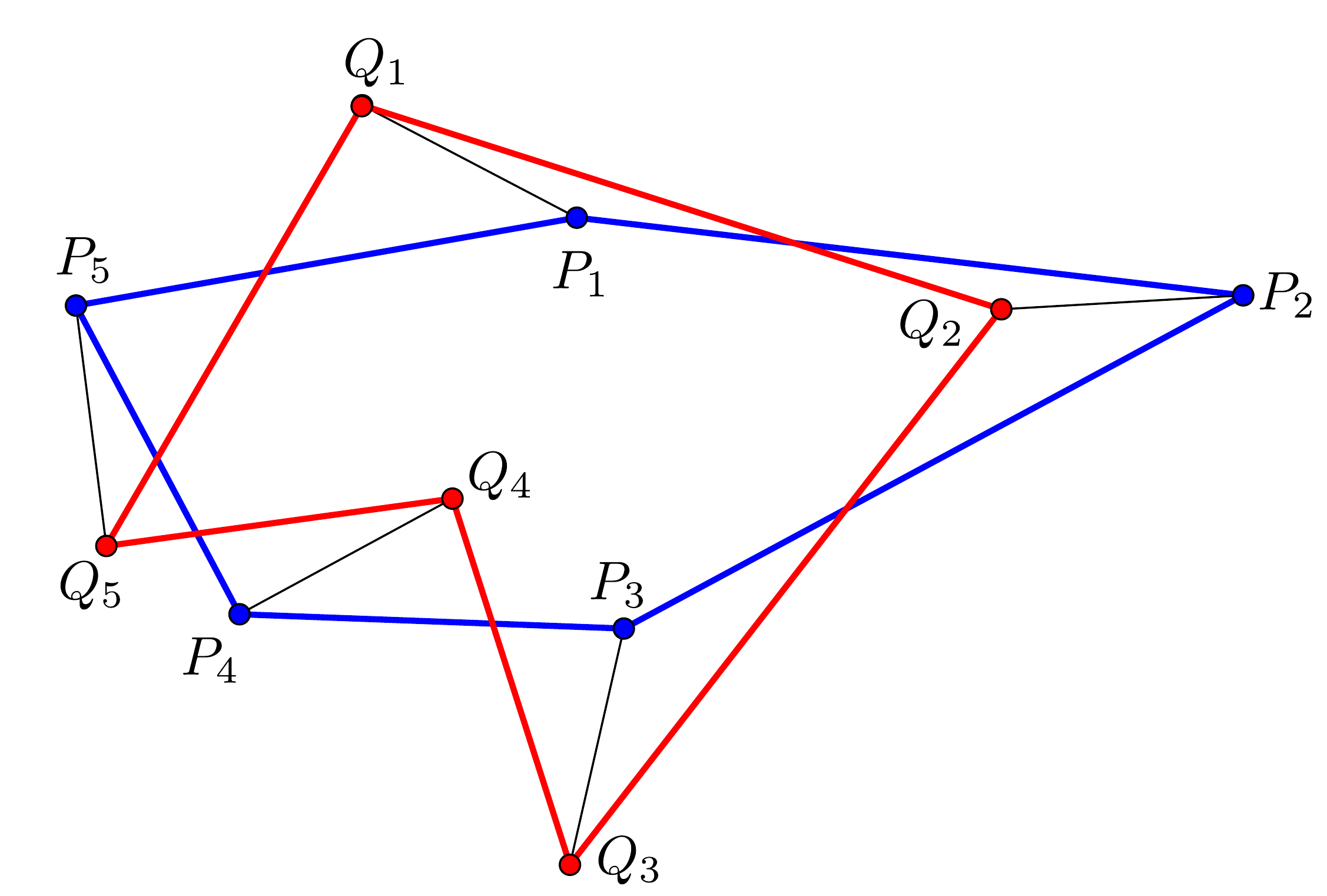}
\caption{Left: folding the parallelogram $P_1Q_1Q_1'P_2$ to the trapezoid $P_1Q_1P_2Q_2$. Right: the pentagons ${\bf P}$ and ${\bf Q}$ are in the discrete bicycle correspondence.}	
\label{chain}
\end{figure}

Two $n$-gons, ${\bf P}=(P_1,\ldots,P_n)$ and ${\bf Q}=(Q_1,\ldots,Q_n)$, are in the discrete bicycle correspondence if every quadrilateral $P_i Q_i  Q_{i+1} P_{i+1}$ is 
obtained by folding a parallelogram along a diagonal as shown in Figure \ref{chain}:
\begin{equation} \label{eq:ldef}
 |P_iP_{i+1}|=|Q_iQ_{i+1}|,\ |P_i,Q_i|=c,\ i=1,\ldots,n, 
\end{equation}
where $c$ is a fixed parameter. In other words, $P_i Q_i  P_{i+1} Q_{i+1}$ is an equilateral trapezoid, perhaps self-intersecting.

The discrete bicycle correspondence is completely integrable. Specifically, in \cite{TT} a Lax presentation with a spectral parameter of the discrete bicycle correspondence is described, providing integrals of this correspondence.

\begin{figure}[ht]
\centering
\includegraphics[width=.45\textwidth]{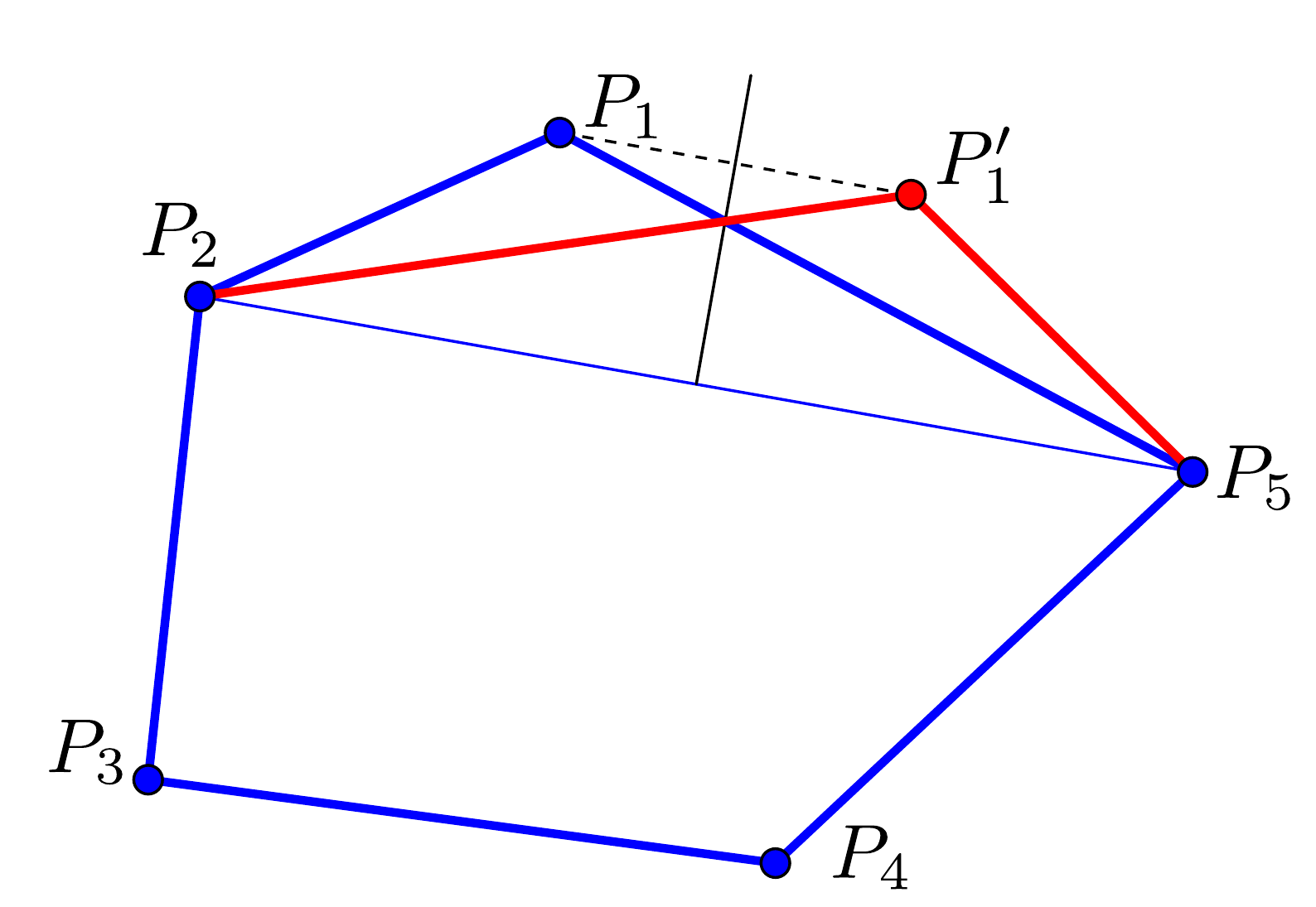}
\caption{Polygon recutting on vertex $P_1$: point $P_1'$ is the reflection of $P_1$ in the perpendicular bisector of the diagonal $P_2P_5$. The recutting of the polygon  is the result of five such transformations performed cyclically.}	
\label{Adler}
\end{figure}

It is also shown there that the discrete bicycle correspondence commutes and shares its integrals with the {\it polygon recutting}, another integrable transformation of polygons, introduced and studied in \cite{Ad1,Ad2}, see Figure \ref{Adler}. 

The present paper concerns an analog of the discrete bicycle correspondence in the centroaffine geometry, associated with the group $SL(2,\R)$ -- or $SL(2,\C)$, if one works with complex coefficients. That is, we consider two polygons in $\R^2$ congruent if they are related by a linear transformation with determinant 1. We denote the determinant by  bracket.

\begin{figure}[ht]
\centering
\includegraphics[width=.55\textwidth]{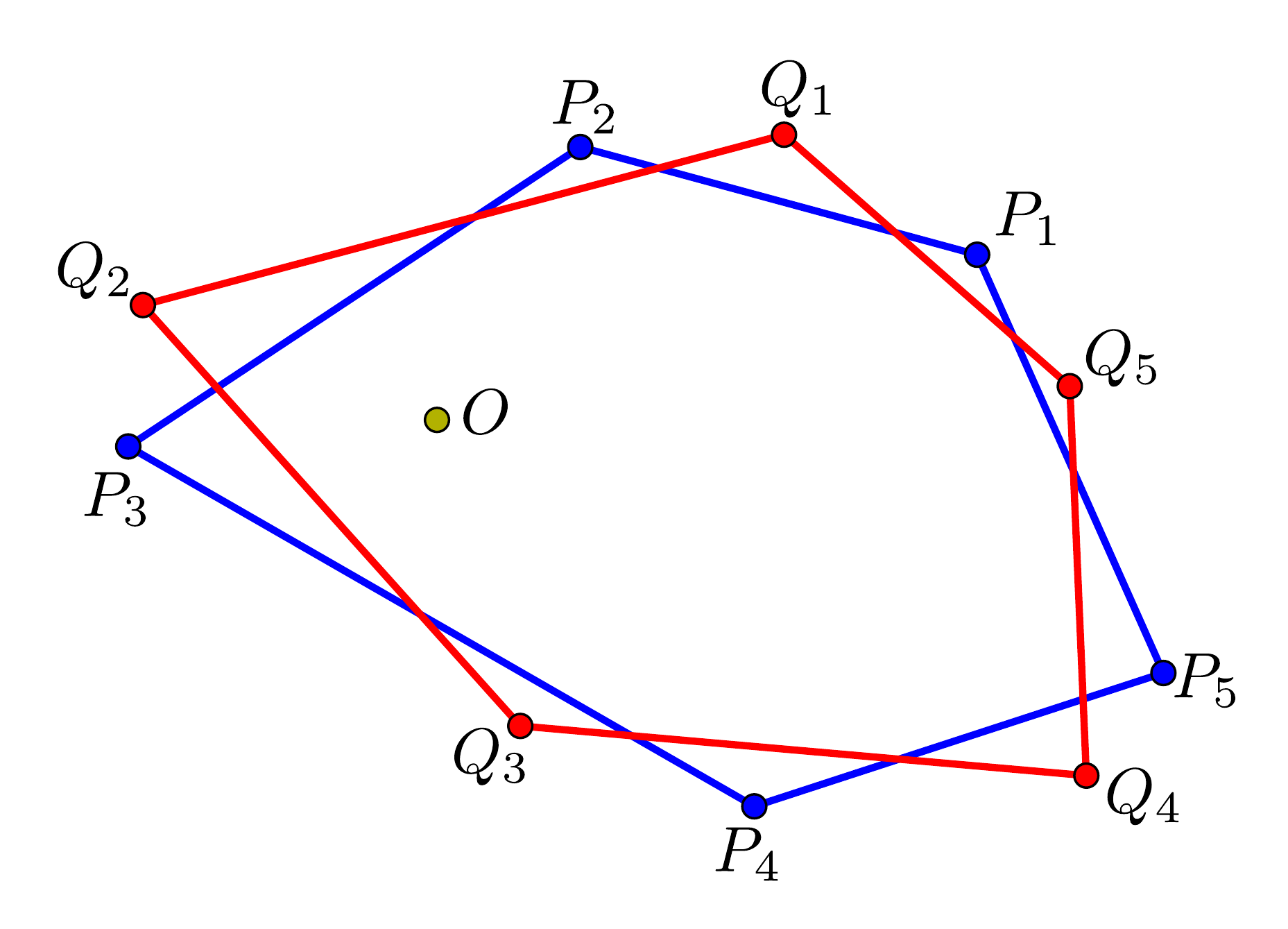}
\caption{Two $c$-related pentagons.}	
\label{pentagons}
\end{figure}

Let $c\in\R$ be a non-zero number.
Two $n$-gons, ${\bf P}$ and ${\bf Q}$, are {\it $c$-related} if 
\begin{equation} \label{eq:cdef}
[P_i, P_{i+1}] = [Q_i, Q_{i+1}],\  [P_i, Q_i]=c,\ i=1,\ldots,n,
\end{equation}
see Figure \ref{pentagons}. We write ${\bf P} \crel {\bf Q}$.

Equations (\ref{eq:cdef}) are  centroaffine analogs of equations (\ref{eq:ldef}): the role of the length is played by the area (i.e., the determinant). The space of centroaffine  polygons is foliated by the $c$-relation invariant subspaces consisting of the polygons whose ``side areas" $[P_i, P_{i+1}]$ depend only on $i$. 

Along with closed polygons ($P_{i+n}=P_i$ for all $i$), we  consider {\it twisted $n$-gons}. A twisted $n$-gon ${\bf P}$ is an infinite collection of points $P_i \in \R^2$ such that $P_{i+n}={\mathcal M}_{\bf P}(P_i)$ for all $i$; this map  ${\mathcal M}_{\bf P} \in SL(2,\R)$ is called the monodromy of the twisted polygon ${\bf P}$. 
Twisted polygons ${\bf P}$ and ${\bf Q}$  are {\it $c$-related} if, in addition to (\ref{eq:cdef}), they share their monodromies. 

The $c$-relation is a discretization of a relation on centroaffine curves, which is a geometrical realization of the B\"acklund transformation of the KdV equation studied in \cite{BBT,Ta18}.

\begin{figure}[ht]
\centering
\includegraphics[width=.5\textwidth]{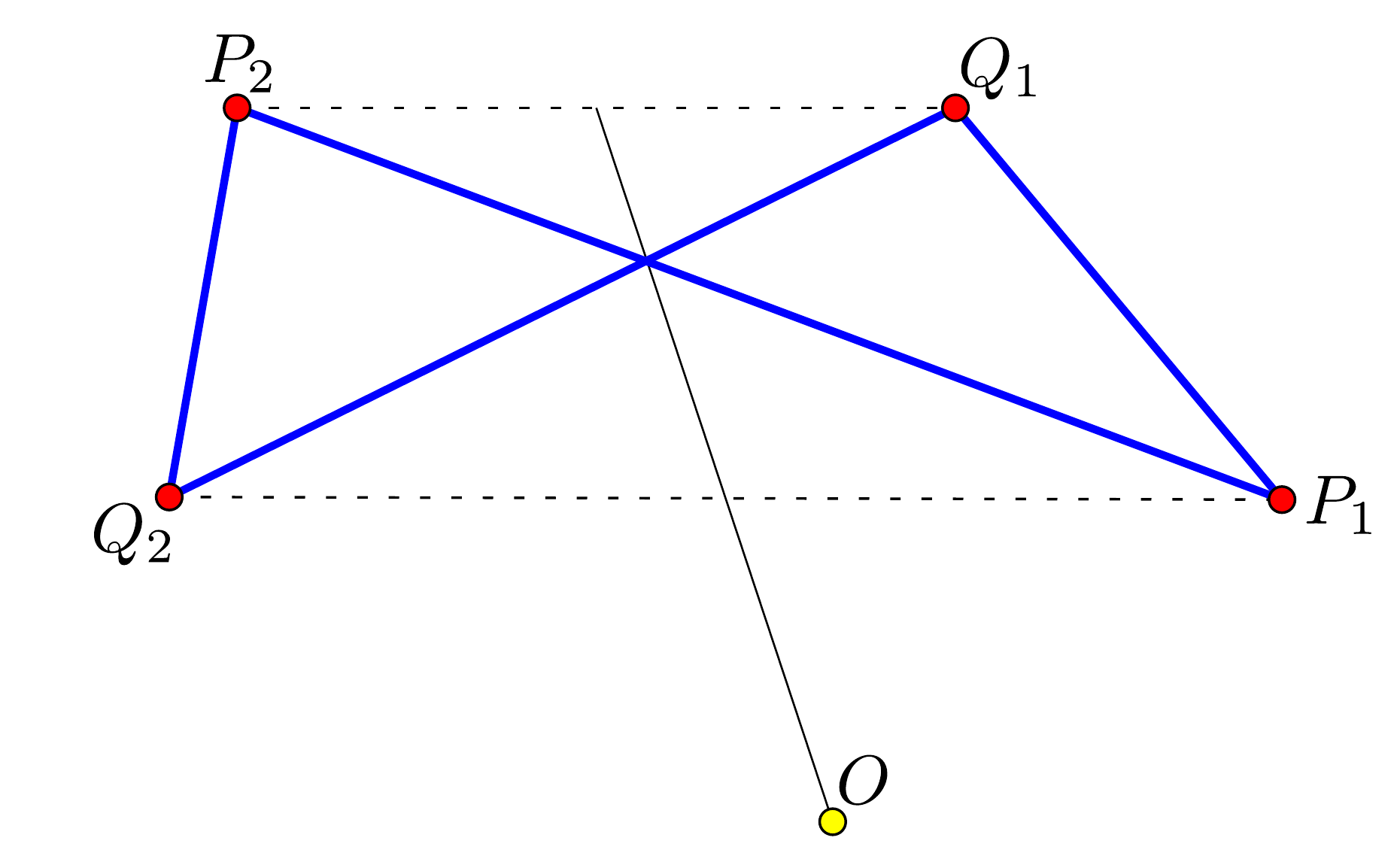}
\caption{A centroaffine butterfly: an affine reflection interchanges  $P_1$ with $Q_2$ and $P_2$ with $Q_1$.}	
\label{butterfly}
\end{figure}

Two consecutive pairs of vertices of $c$-related polygons form a quadrilateral satisfying 
$$
[P_i,P_{i+1}]=[Q_i,Q_{i+1}],\ [P_i,Q_i]=[P_{i+1},Q_{i+1}].
$$
 We call such quadrilaterals {\it centroaffine butterflies}, see Figure \ref{butterfly} (the term is adopted from \cite{TT}). They are centroaffine analogs of the folded parallelograms in Figure  \ref{chain}. 
 
 We also consider centroaffine version of polygon recutting. An elementary centroaffine recutting is depicted in Figure \ref{butterfly}: it is a linear involution that swaps the triangles $P_1Q_1Q_2$ and $P_1P_2Q_2$. The centroaffine recutting of a $n$-gon is  the composition of $n$ elementary recuttings performed cyclically. 
 
 Let $\pi: \R^2\nobreak\setminus\nobreak\{O\} \to \RP^1$ be the natural projection. Abusing notation, we use the same symbol for the projections of centroaffine polygons to polygons in $\RP^1$. This projection commutes with the natural actions of $SL(2,\R)$ on $\R^2$ and $\RP^1$. 
 
 Let ${\bf P}$ and ${\bf Q}$ be $c$-related centroaffine polygons, and ${\bf p}$ and ${\bf q}$ be their projections to $\RP^1$. Then 
\begin{equation} \label{eq:rel}
[p_i,p_{i+1},q_i,q_{i+1}]=\frac{[P_i,P_{i+1}][Q_i,Q_{i+1}]}{[P_i,Q_i][P_{i+1},Q_{i+1}]}
\end{equation}
where the bracket on the left hand side denotes the cross-ratio (there are six different choices of cross-ratio to make; the right hand side of the formula specifies our choice).
If $[P_i,P_{i+1}]$ is the same for all $i$, then the polygons ${\bf p}$ and ${\bf q}$ are in the cross-ratio relation: $[p_i,p_{i+1},q_i,q_{i+1}]=\alpha$ for all $i$. 

The cross-ratio relation on projective polygons was thoroughly studied, starting with \cite{NC} and, more recently, in \cite{AFIT} and  \cite{AGR}. This is the second source of our motivation: many results in this paper have analogs in \cite{AFIT}. 

The cross-ratio relation can be generalized: $n$-gons ${\bf p}$ and ${\bf q}$ are related if $[p_i,p_{i+1},q_i,q_{i+1}]=\alpha_i$, where $\alpha_i$  is an $n$-periodic sequence (not necessarily constant). Formulas (\ref{eq:cdef}) and  (\ref{eq:rel}) imply that the projection $\pi$ conjugates the $c$-relation with this generalized cross-ratio relation.
 \medskip
 
Let us present the main results of the paper. 

In Section \ref{sec:spaces} we introduce coordinates in the moduli space of twisted centroaffine polygons and calculate the monodromy of a twisted polygon: the result is given in terms of continuants (3-diagonal determinants). We also describe the algebraic relations satisfied by the coordinates of closed polygons.

Let ${\bf P}=(P_1,\ldots,P_n)$ be a centroaffine  $n$-gon. Choose a test vector $Q_1$ with $[P_1,Q_1]=c$, and consecutively construct vectors $Q_2,\ldots,Q_n,Q_{n+1}$ according to (\ref{eq:cdef}). We call the map $Q_1\mapsto Q_{n+1}$ the  {\it Lax transformation} associated with the polygon ${\bf P}$ and denote it by ${\mathcal L}_{{\bf P},c}$. Then ${\bf P} \crel {\bf Q}$ if and only if $Q_1$ is a fixed point of ${\mathcal L}_{{\bf P},c}$. Similarly one defines the Lax transformation associated with a twisted polygon.
 
 We show that ${\mathcal L}_{{\bf P},c}$ is a M\"obius map. This makes it possible to consider the $c$-relation on generic polygons as a 2-2 map.
 
 Theorem \ref{thm:Lax} states that if 
${\bf P}$ and ${\bf Q}$ are $c$-related twisted $n$-gons, then the Lax transformations ${\mathcal L}_{{\bf P},\lambda}$ and ${\mathcal L}_{{\bf Q},\lambda}$ are conjugated for every value of the spectral parameter $\lambda$. This is the source of integrals of the $c$-relation. The moduli space of twisted centroaffine $n$-gons with fixed ``side areas" has dimension $n$; we obtain $\lfloor{\frac{n+1}{2}}\rfloor$ integrals therein.

Theorem \ref{thm:Bianchi} states that the $c$-relation satisfies the Bianchi permutability. Informally speaking, it says that the $c$-relations with different values of the constant $c$ commute (see Section \ref{subsec:Bianchi} for the precise formulation).

 In Section \ref{subsec:recut} we study how the $c$-relation interacts with the centroaffine polygon recutting. We prove that 
 the Lax transformation is preserved by the recutting and that the recutting commutes with the $c$-relations (Theorem \ref{thm:recut}). 	
 
 In Section \ref{subsec:inttwist} we calculate the integrals provided by the conjugacy invariance of the Lax transformations and show that they coincide with the integrals of the dressing chain of Veselov and Shabat \cite{VS}. In Section \ref{subsec:intclosed} we describe two relations between these integrals that hold for closed polygons.
 
Section \ref{sec:before} concerns the space of closed centroaffine polygons before its factorization by the group $SL(2,\R)$. We construct presymplectic forms on the subspaces of polygons whose ``side areas" $[P_i,P_{i+1}]$ depend on $i$ only. These forms are $SL(2,\R)$-invariant, but they do not descend on the quotient spaces by the group. The forms are invariant under the $c$-relation and under the polygon recutting.  

We  construct three additional integrals of the $c$-relation, quadratic in the coordinates; a polynomial function of these three integrals is invariant under $SL(2,\R)$. These integrals are interpreted as the moment map of the Hamiltonian action of $sl(2,\R)$ on  the spaces of polygons with fixed ``side areas". 

These integrals define a certain {\it center} $C({\bf P})$ of a centroaffine polygon, that takes values in quadratic forms on $\R^2$.
This center is  invariant under the $c$-relation and the recutting, and equivariant with respect to the action of $SL(2,\R)$.
It has interesting properties: it is additive with respect to cutting polygons into two, and it coincides with the origin for centroaffine butterflies. 

Section \ref{sec:small} is devoted to a study of ``small-gons", triangles, quadrilaterals, and pentagons. 
 
We emphasize that other aspects of complete integrability, such as invariant Poisson structures and relation with the theory of cluster algebras, are not discussed in the present paper. They will be studied by A. Izosimov in the forthcoming paper \cite{Iz}.

\bigskip

{\bf Acknowledgements}. We are very grateful to Anton Izosimov for his insights and useful suggestions. 
 ST was supported by NSF grant DMS-2005444. 
 

\section{Spaces and maps} \label{sec:spaces}

\subsection{Spaces and coordinates} \label{subsec:spaces}

In this paper we consider polygons ${\bf P}=(\ldots P_i P_{i+1} \ldots)$ in $\R^2$ that satisfy $[P_i,P_{i+1}]\neq 0$ for all $i$ (when appropriate, the indices are understood cyclically). 
Denote by $\widetilde {\mathcal Y}_n$ and $\widetilde {\mathcal X}_n$ the spaces of twisted and closed  $n$-gons, and by ${\mathcal Y}_n$ and ${\mathcal X}_n$ their quotient spaces by $SL(2,\R)$. 

 Let us introduce  coordinates in ${\mathcal Y}_n$:
\begin{equation}
	\label{eq:coordinates}
	\ed_{2j-1}=[P_{j-1},P_j],\ \ve_{2j}=[P_{j-1},P_{j+1}].
\end{equation}
That is, $\ed_{2j-1}$ are the areas subtended by the sides, and $\ve_{2j}$ are the areas subtended by the short diagonals of the polygon.

One has a linear recursion
\begin{equation} \label{eq:P_j+1}
P_{j+1}= \frac{\ve_{2j}}{\ed_{2j-1}}P_j - \frac{\ed_{2j+1}}{\ed_{2j-1}}P_{j-1},
\end{equation}
that is,
$$
\begin{pmatrix}
P_{j}\\
P_{j+1} 
\end{pmatrix}
=
\begin{pmatrix}
0 & 1\\
- \frac{\ed_{2j+1}}{\ed_{2j-1}} & \frac{\ve_{2j}}{\ed_{2j-1}}\\
\end{pmatrix}
\begin{pmatrix}
P_{j-1}\\
P_{j} 
\end{pmatrix}.
$$
It follows that the (conjugacy class of the) monodromy is given by
\begin{equation} \label{eq:monod}
{\mathcal M}_{\bf P} = 
\begin{pmatrix}
0 & 1\\
 - \frac{\ed_{1}}{\ed_{2n-1}} & \frac{\ve_{2n}}{\ed_{2n-1}} \\
\end{pmatrix}
\ldots
\begin{pmatrix}
0 & 1\\
- \frac{\ed_{5}}{\ed_{3}} & \frac{\ve_4}{\ed_{3}}\\
\end{pmatrix}
\begin{pmatrix}
0 & 1\\
- \frac{\ed_3}{\ed_{1}} & \frac{\ve_2}{\ed_{1}}\\
\end{pmatrix}.
\end{equation}
Note that $\det {\mathcal M}_{\bf P} =1$.

Let ${\bf S}=(s_1,s_3,s_5,\ldots, s_{2n-1})$; denote by $\widetilde{\mathcal Y}_{n,{\bf S}}$ the space of twisted $n$-gons with $[P_i,P_{i+1}]=s_{2i+1}$ for all $i$, and $\widetilde{\mathcal X}_{n,{\bf S}}$ for closed polygons. 

\begin{remark}
{\rm The spaces $\widetilde{\mathcal X}_{n,{\bf S}}$ and ${\mathcal X}_{n,{\bf S}}$ are centroaffine analogs of the spaces of polygons with fixed side lengths, studied in \cite{KM}.
}
\end{remark}

When is $\widetilde{\mathcal X}_{n,{\bf S}}$ a smooth $n$-dimensional manifold? One has a map $\widetilde{\mathcal X}_{n} \to \R^n$ that sends ${\bf P}$ to ${\bf S}=(\ldots,[P_i,P_{i+1}],\ldots)$. The next lemma describes the regular values ${\bf S}$ of this map.

\begin{lemma} \label{lm:regval}
If $n$ is odd and $s_j\ne 0$ for all $j$, then  $\widetilde{\mathcal X}_{n,{\bf S}}$ is smooth. If $n$ is even and $s_1 s_5 s_9 \cdots \neq \pm s_3 s_7 s_{11} \cdots$, then  $\widetilde{\mathcal X}_{n,{\bf S}}$ is smooth. 
\end{lemma}

\begin{proof}
If $P_i=(x_i,y_i)$ then $s_{2i+1}=x_iy_{i+1}-x_{i+1}y_i$. We need to know when the 1-forms $ds_{2i+1}$ are linearly dependent. 

We may assume that the coordinates are chosen so that all $x_i$ and $y_i$ are distinct from zero. Assume that
$$
\sum_i a_i (x_i dy_{i+1} + y_{i+1}dx_i-x_{i+1}dy_i-y_idx_{i+1}) =0,
$$
where not all $a_i$ vanish, hence
$$
\sum_i (a_i y_{i+1} - a_{i-1}y_{i-1}) dx_i - (a_i x_{i+1} - a_{i-1}x_{i-1}) dy_i =0.
$$
Therefore
$$
\frac{a_{i-1}}{a_i}=\frac{y_{i+1}}{y_{i-1}}=\frac{x_{i+1}}{x_{i-1}},\ i=1,\ldots,n,
$$
(and hence all $a_i$ are different from zero).

This implies that $P_{i-1}$ and $P_{i+1}$ are collinear for all $i$. If $n$ is odd, then all vertices of ${\bf P}$ are collinear, and $s_{2i+1}=0$ for all $i$.

If $n$ is even, then the odd vertices of ${\bf P}$ are collinear, and so are the even vertices. Without loss of generality, assume that
$$
P_1=(x_1,0), P_2=(0,y_2), P_3=(x_3,0), P_4=(0,y_4),\ldots
$$ 
Then 
$$
s_3=x_1y_2, s_5=-y_2x_3, s_7=x_3y_4, s_9=-y_4x_5 \ldots
$$
and $s_3 s_7 \cdots = \pm s_1 s_5 \cdots$
This completes the proof.
\end{proof}

\subsection{Monodromy of a twisted polygon} \label{subsec:mono}

Set
$$
a_j=\frac{v_{2j}}{s_{2j-1}},\ b_j=\frac{s_{2j+1}}{s_{2j-1}}.
$$
Then 
$$
{\mathcal M}_n = 
\begin{pmatrix}
0&1\\
 - b_{n-1}&a_{n-1} 
\end{pmatrix}
\ldots
\begin{pmatrix}
0&1\\
- b_1&a_1
\end{pmatrix}
\begin{pmatrix}
0&1\\
- b_0&a_0
\end{pmatrix},
$$
as in (\ref{eq:monod}). 

For $0\le i\le j$, consider the  continuants
\begin{equation}
	\label{eq:continuant_ij}
D_{i,j+1}=	
\det     
\begin{pmatrix}
	a_i&1&0&\cdots&0\\
	b_{i+1}&a_{i+1}&1&0&\cdots\\
	0&b_{i+2}&a_{i+2}&1&\cdots\\
	0&\ddots&\ddots&\ddots&0\\
	\cdots&0&b_{j-1}&a_{j-1}&1\\
	0&\cdots&0&b_j&a_j
\end{pmatrix}.
\end{equation}
We describe the monodromy of a twisted polygon in terms of these continuants.

Expanding by the last two rows, we see that
\begin{equation} \label{eq:recD}
D_{i,j+1}=   a_jD_{i,j} - b_jD_{i,j-1},
\end{equation}
the same recursion as (\ref{eq:P_j+1}).
Let us add the boundary conditions
$
D_{i,i}=1,\ D_{i,i-1}=0. 
$

\begin{proposition} \label{prop:monod}
One has
$$
{\mathcal M}_n = 
\begin{pmatrix}
-b_0 D_{1,n-1} & D_{0,n-1}\\
-b_0 D_{1,n} & D_{0,n}
\end{pmatrix}.
$$
\end{proposition}

\begin{proof}
Induction on $n$. For $n=1$, the claim holds due to the boundary conditions
$D_{1,0}=0, D_{0,0}=D_{1,1}=1, D_{0,1}=a_0$.

Next, 
$$
{\mathcal M}_{n+1} =
\begin{pmatrix}
0&1\\
 - b_{n}&a_{n} 
\end{pmatrix}
{\mathcal M}_{n},
$$
and the result follows from the recurrence (\ref{eq:recD}). 
\end{proof}

If ${\bf P}$ is a twisted $n$-gon, then the sequences $a_i$ and $b_i$ are $n$-periodic and $\prod_i b_i =1$.
Proposition \ref{prop:monod} implies the following statement.

\begin{corollary} \label{cor:closed}
For a closed $n$-gon, the coordinates $(s_{2i-1},v_{2i})$ satisfy
$D_{i,n+i-1} = 0$ for all $i$.
\end{corollary}

In fact, any three of these identities imply the rest (the codimension of the space of closed polygons is three).

\begin{example} \label{ex:34}
{\rm Consider the case $n=3$. Corollary \ref{cor:closed} implies
$$
\ve_0\ve_2=\ed_3\ed_5,\ \ve_2\ve_4=\ed_5\ed_1,\ \ve_4\ve_0=\ed_1\ed_3,
$$
hence $\ve_0^2= \ed_3^2$. We have two solutions;
$$
\ve_0=-\ed_3,\  \ve_2=-\ed_5,\ \ve_4=-\ed_1,\ {\rm and}\ \ \ve_0=\ed_3,\  \ve_2=\ed_5,\ \ve_4=\ed_1.
$$
The first one corresponds to a closed triangle (the monodromy is $Id$),  and the second one to a centrally symmetric hexagon (the monodromy is $-Id$).

Next, consider the case $n=4$. Corollary \ref{cor:closed} implies
$$
\ve_0\ve_2\ve_4=\ed_3\ed_7\ve_4+\ed_1\ed_5\ve_0
$$
and its three cyclic permutations. Rewrite it as
$$
\frac{\ed_3\ed_7}{\ve_0\ve_2} + \frac{\ed_1\ed_5}{\ve_2\ve_4}=1
$$
and its cyclic permutations. This is a system of four linear equations on the variables
$$
\frac{1}{\ve_0\ve_2},\ \frac{1}{\ve_2\ve_4},\ \frac{1}{\ve_4\ve_6},\ \frac{1}{\ve_6\ve_0}
$$
with coefficients $\ed_3\ed_7$ and $\ed_1\ed_5$. This system implies
$$
\frac{1}{\ve_0\ve_2} = \frac{1}{\ve_4\ve_6},\ \frac{1}{\ve_2\ve_4} = \frac{1}{\ve_6\ve_0},
$$
and hence $\ve_4^2=\ve_0^2$. As before, one has two choices of signs, one corresponding to the monodromy  $Id$, and another to $-Id$.
In the former case of closed quadrilaterals, one has 
$$
\ve_4=-\ve_0, \ve_6=-\ve_2\ {\rm and}\ \ \ve_0\ve_2=\ed_7\ed_3-\ed_1\ed_5,
$$
the latter being the Ptolemy-Pl\"ucker relation.
}
\end{example}

\subsection{Lax transformation is fractional-linear
}  \label{subsec:maps}

Given a non-zero vector $P_i$, the $i$th vertex of a polygon ${\mathbf P}$, 
the vectors $Q$ with $[P_i,Q]=c$ comprise a line $L_{P_i}$ parallel to $P_i$. Identify $L_{P_i}$ with $\R$ by  parameterizing it as 
$$
\frac{cP_{i+1}}{[P_i,P_{i+1}]} + tP_i, \ t\in\R.
$$

Let ${\bf P}=(P_1,\ldots,P_n)$ be a closed  $n$-gon. As described in Introduction, 
choose a test vector $Q_1$ with $[P_1,Q_1]=c$, and consecutively construct vectors $Q_2,\ldots,Q_n,Q_{n+1}$ according to (\ref{eq:cdef}). We  use the notation 
$$
{\mathcal L}_{P_1 P_2,c} (Q_1)=Q_2,\ {\mathcal L}_{P_2 P_3,c} (Q_2)=Q_3, \ldots
$$
(we omit $c$ from the notation when it does not lead to confusion).
The   Lax transformation  ${\mathcal L}_{{\bf P},c}: L_{P_1}\to L_{P_1}$ is the composition of these maps. Then ${\bf P} \crel {\bf Q}$ if and only if $Q_1$ is a fixed point of ${\mathcal L}_{{\bf P},c}:\R\to \R$. 
 
 Likewise, if ${\mathbf P}$ is a twisted $n$-gon with monodromy ${\mathcal M}_{\bf P}$, then  ${\mathcal M}_{\bf P}$ sends $L_{P_1}$ to $L_{P_{n+1}}$ and, as a map $\R\to\R$, it is the identity. Hence it is still true that the fixed points of ${\mathcal L}_{{\bf P},c}$ give rise to twisted polygons ${\mathbf Q}$ such that ${\bf P} \crel {\bf Q}$.

\begin{lemma} \label{lm:Lax}
The Lax transformations are fractional-linear.
\end{lemma}

\begin{proof}
Denote by $R_{Q_1,P_2}$ the centroaffine reflection that interchanges $Q_1$ and $P_2$.
 This map is given by the formula
\begin{equation} \label{eq:invol}
R_{Q_1,P_2} (X)=\frac{[Q_1,X]Q_1+[X,P_2]P_2}{[Q_1,P_2]},
\end{equation}
and one has $Q_2=R_{Q_1,P_2} (P_1)$, that is,  ${\mathcal L}_{P_1 P_2} (Q_1) = R_{Q_1,P_2} (P_1)$.

Using the identifications of $L_{P_1}$ and $L_{P_2}$ with $\R$, the map ${\mathcal L}_{P_1 P_2}: Q_1 \mapsto Q_2$ becomes
$$
t \mapsto -\frac{c[P_1,P_3]}{[P_1,P_2][P_2,P_3]} + \frac{[P_1,P_2]^2-c^2}{t[P_1,P_2]^2} 
= -\frac{c\ve_4}{\ed_3\ed_5} + \frac{\ed_3^2-c^2}{t\ed_3^2},
$$
a fractional-linear transformation.
Hence ${\mathcal L}_{{\bf P},c}$ is fractional-linear as well.
\end{proof}

The above fractional-linear transformation is represented by the matrix
\begin{equation} \label{eq:Laxmat}
\begin{pmatrix}
	-\dfrac{c\ve_4}{\ed_{3}\ed_{5}}&\displaystyle{1-\frac{c^2}{\ed_3^2}}\\
	1&0
\end{pmatrix},
\end{equation}
 its determinant equals $\displaystyle{\frac{c^2}{\ed_{3}^2}-1}$.
 
 If the ground field is $\C$, a fractional-linear transformation has two fixed points, perhaps, coinciding (unless it is the identity). Over the reals, this number is  0, or 1, 2, or $\infty$.

Over $\C$, Lemma \ref{lm:Lax} makes it possible to consider $c$-relation as a map (defined in a Zariski open set of polygons): start with ${\mathbf P}$ and choose one of the two $c$-related polygons, say ${\mathbf Q}_1$. This polygon also has two $c$-related ones, one of which is $-{\mathbf P}$, the polygon that is centrally symmetric to ${\mathbf P}$. Choose the other polygon $c$-related to ${\mathbf Q}_1$ and continue in a similar way. If at the beginning we choose ${\mathbf Q}_2$ instead of ${\mathbf Q}_1$ then $-{\mathbf Q}_2 \crel {\mathbf P}$, that is, up to central symmetry, we obtain the inverse of the same map. 

Over $\R$, a polygon may have no $c$-related ones. To have a map, one needs to assume that a $c$-related polygon exists.
We shall see in Section \ref{subsec:inv} that the $c$-related polygons have conjugated Lax transformations. This makes it possible 
to continue in the same way as in the complex case. 
 
\subsection{Centroaffine butterflies} \label{subsec:but}

\begin{figure}[ht]
\centering
\includegraphics[width=.5\textwidth]{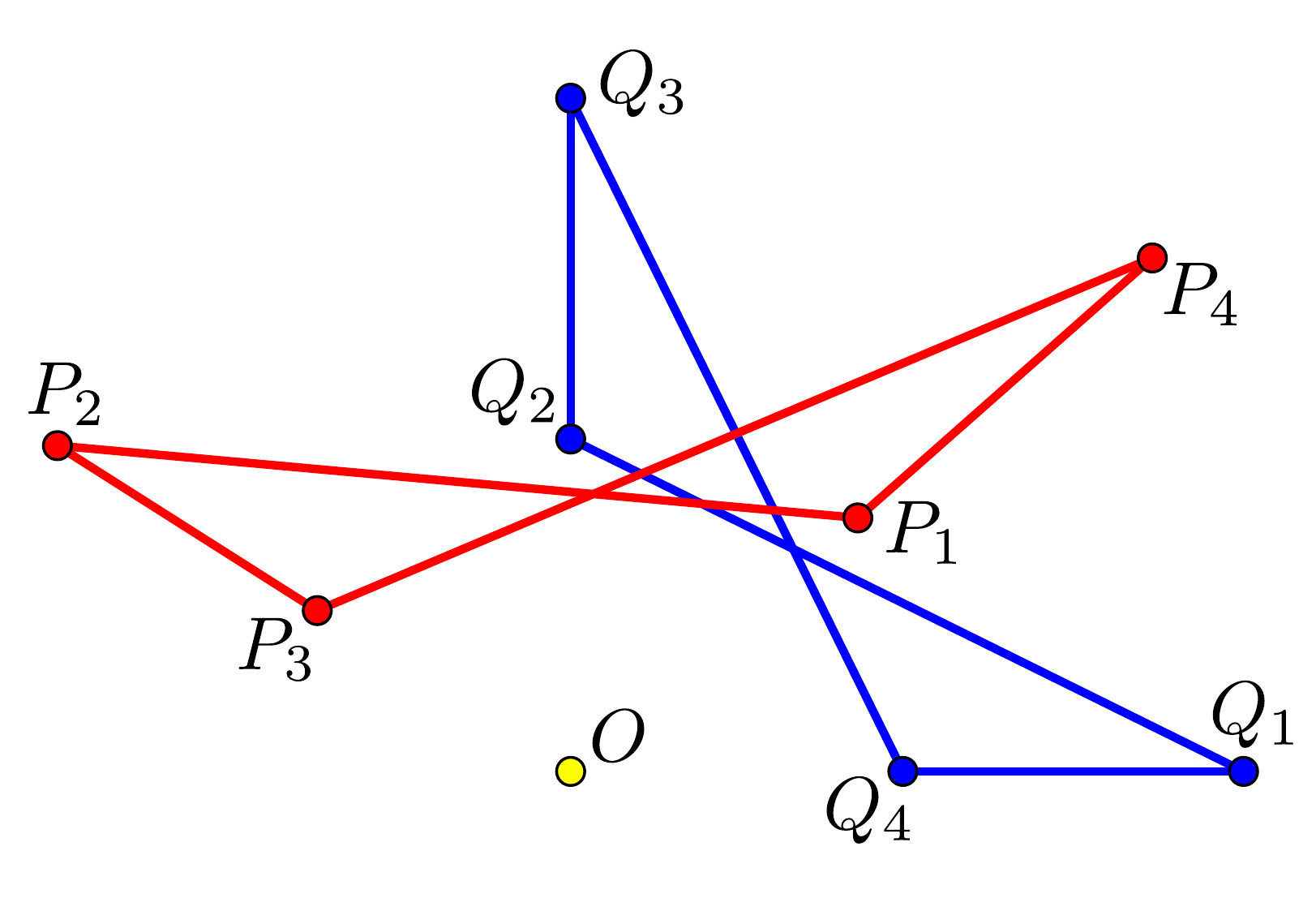}
\caption{Two centroaffine butterflies are $c$-related.}
\label{Butterflies}
\end{figure}

We show that centroaffine butterflies have trivial Lax transformations for all $c$ (see Figure \ref{Butterflies}) and classify all such quadrilaterals. 

\begin{lemma} \label{lm:butter}
The Lax transformation for a quadrilateral is the identity for every $c$ if and only if the quadrilateral is a centroaffine butterfly, or it is obtained from a centroaffine butterfly by reflecting one of the vertices in the origin (an ``anti-butterfly"), or its two opposite vertices are symmetric with respect to the origin.
\end{lemma}

\begin{proof}  Assume first that two opposite vertices of a quadrilateral ${\mathbf P}$ are not collinear. 
 Applying a linear transformation, assume that
$$
P_1=(1,0),\ P_2=(a,b),\ P_3=(0,1),\ P_4=(u,v), \ Q_1=(x,y),
$$
so $y=c$. Applying equation (\ref{eq:invol}) twice,  we find
$$
Q_2=
 \frac{1}{b x-a y}(a b-x y , b^2-y^2 ), \
 Q_3=\left(-y,\frac{-y \left(a^2+b^2\right)+a b x+y^3}{a b-x y}\right).
 $$
 Going in the opposite direction, that is, replacing $(a,b)$ with $(u,v)$, we obtain point
 $$
 \left(-y,\frac{-y \left(u^2+v^2\right)+uv x+y^3}{uv-x y}\right).
 $$ 
 This point coincides with $Q_3$ for all $x$, which implies 
 $$
 ab=uv,\ a^2+b^2=u^2+v^2.
 $$
 Hence either $u=b,v=a$, a butterfly, or $u=-b,v=-a$, an anti-butterfly, or $u=-a,v=-b$, symmetric opposite vertices. 
 
 It remains to consider the case when both pairs of opposite vertices of ${\mathbf P}$ are collinear. Then one may assume that
 $$
 P_1=(1,0),\ P_2=(0,1),\ P_3=(a,0),\ P_4=(0,b), \ Q_1=(x,y).
 $$
Using equation (\ref{eq:invol}) consecutively, one calculates 
$$
Q_5=\left(\frac{x(a^2-y^2)(b^2-y^2)}{(a^2b^2-y^2)(1-y^2)},y\right).
$$
Equating it to $Q_1$ yields 
$$
(a^2-y^2)(b^2-y^2)=(a^2b^2-y^2)(1-y^2),
$$
hence $(a^2-1)(b^2-1)=0.$ If $a=1$ or  $b=1$, then $P_1=P_3$ or $P_2=P_4$, and if 
$a=-1$ or $b=-1$, we have the already considered symmetric case. 
\end{proof}

\subsection{Constructing $c$-related  polygons} \label{subsec:constr}

In this section we describe a construction that yields a pair of $c$-related $n$-gons. This is a centroaffine analog of a construction that yields pairs of polygons in the discrete bicycle correspondence that is described in \cite{TT}.

\begin{figure}[ht]
\centering
\includegraphics[width=.7\textwidth]{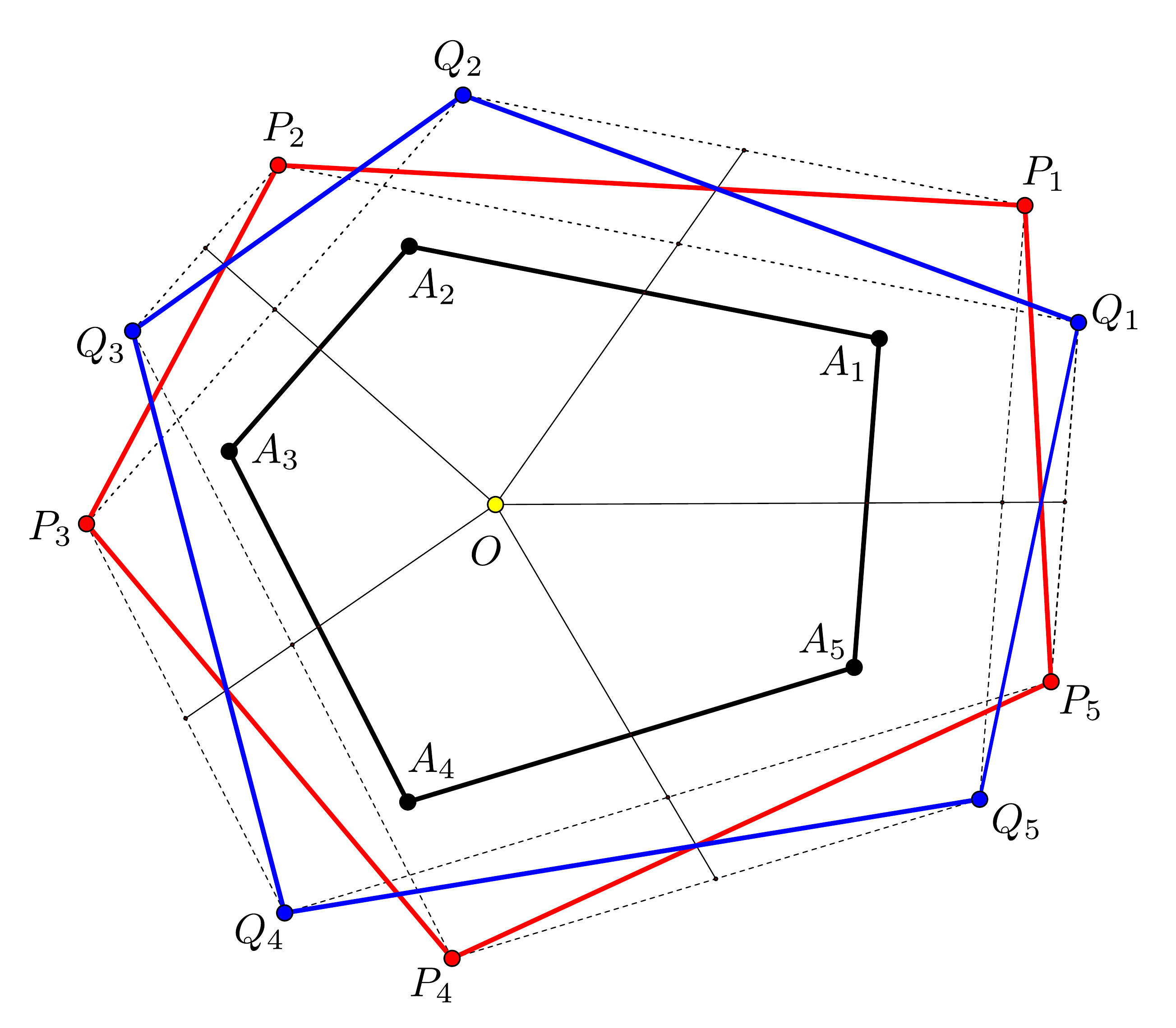}
\caption{Constructing a $c$-related pair of $n$-gons, odd $n$.}
\label{reflections}
\end{figure}

Assume first that $n$ is odd. Start with an $n$-gon (pentagon ${\mathbf A}$ in Figure \ref{reflections}). Connect the midpoints of its sides with the origin. Consider the affine reflections in these lines that interchange the vertices of the respective sides of ${\mathbf A}$. Let $R$ be the composition of these reflections taken around the polygon.

\begin{lemma} \label{lm:reflodd}
The map $R$ is an affine reflection. 
\end{lemma}

\begin{proof}
Since $n$ is odd, $R$ is orientation reversing and its determinant equals $-1$. Hence it has two eigendirections. Also $R$ has a fixed point, a vertex of  ${\mathbf A}$. Therefore the eigenvalues of $R$ are $1$ and $-1$, and it is an affine reflection.
\end{proof} 

It follows that $R^2 = Id$. The construction of a $c$-related pair follows.

Start with an arbitrary point $P_1$ and apply consecutive affine reflections around the polygon ${\mathbf A}$ twice. This produces $2n$ points
$$
P_1 \mapsto Q_2 \mapsto P_3 \mapsto Q_4 \mapsto \cdots \mapsto Q_n,
$$
see \ref{reflections}. Each quadrilateral $P_i P_{i+1} Q_i Q_{i+1}$ is a centroaffine butterfly, therefore ${\bf P} \crel {\bf Q}$. 

One has $c=[P_1,Q_1]=[P_1,R(P_1)]$. The locus of points $P_1$ for which $[P_1,R(P_1)]$ is fixed is a hyperbola (indeed, $R$ is conjugated to the map $(x,y)\mapsto (-x,y)$, in which case this claim is obvious).

Now consider the case of even $n$. We repeat the above construction, but this time the transformation $R$ had determinant 1. It still has an eigendirection with eigenvalue 1 (a vertex of ${\bf A}$ is a fixed point), but it is not necessarily the identity. We need to assume that $R=Id$, see below.


\begin{figure}[ht]
\centering
\includegraphics[width=.8\textwidth]{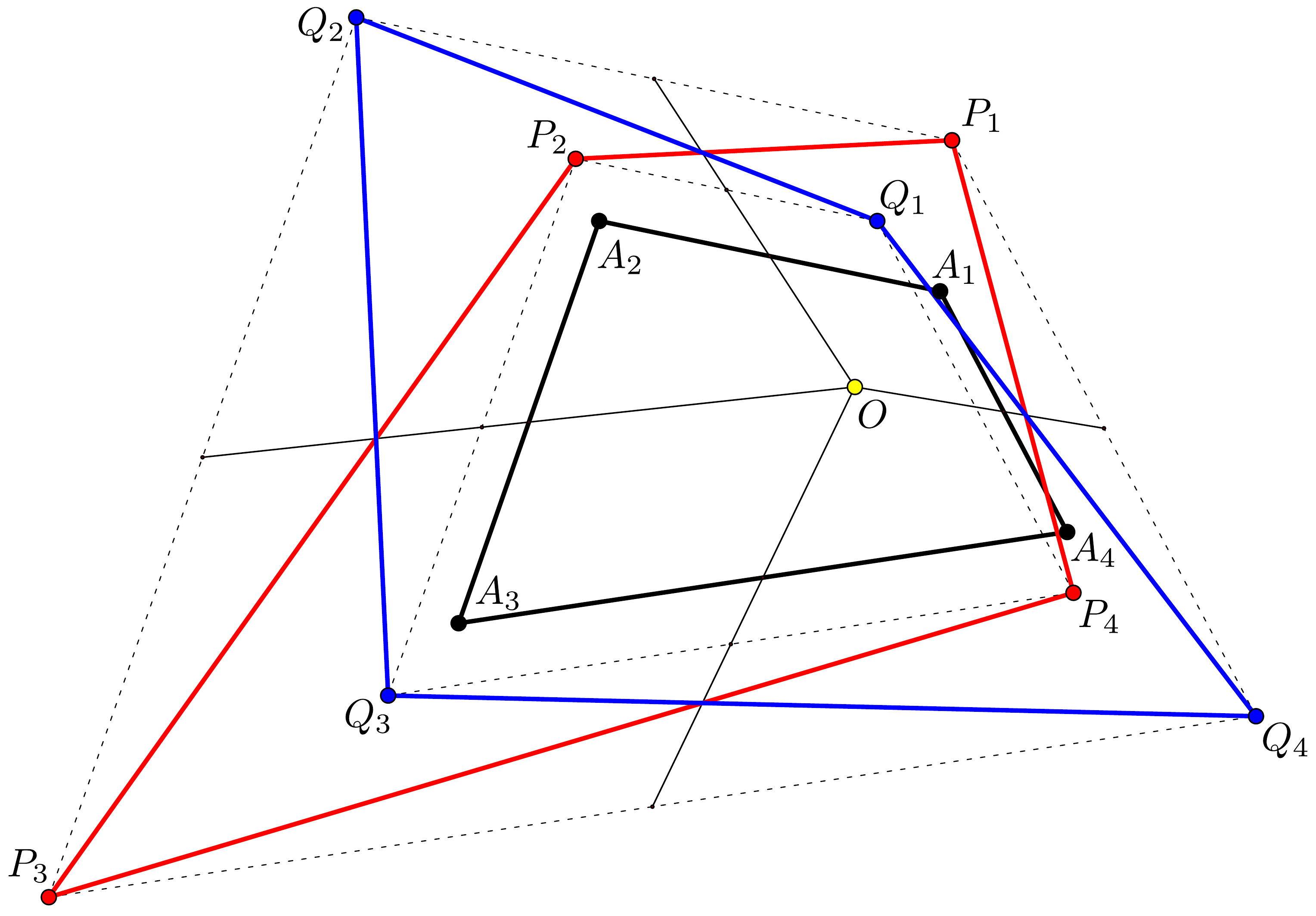}
\caption{Constructing a $c$-related pair of $n$-gons, even $n$.}
\label{reflections2}
\end{figure}

With this assumption, we choose two starting points, $P_1$ and $Q_1$ and apply consecutive affine reflections to obtain  polygons ${\bf P} \crel {\bf Q}$ with $c=[P_1,Q_1]$ as before, see Figure \ref{reflections2}.

We describe when the transformation $R$ is the identity. Recall that we consider the case of an even-gon.

\begin{lemma} \label{lm:ident}
One has $R=Id$ if and only if
$$
\sum_{i=1}^n (-1)^i \frac{\ve_{2i}}{\ed_{2i-1}\ed_{2i+1}} = 0.
$$
\end{lemma} 

\begin{proof}
Since $R$ has a fixed point, a vertex of polygon ${\bf A}$, the transformation $R$ is the identity if and only if it has another fixed point not collinear with the first one. Thus we need to learn when a polygon ${\bf B}$ exists, not homothetic to ${\bf A}$, whose sides are parallel to those of ${\bf A}$ and the midpoints of whose sides are collinear to those of ${\bf A}$.

These conditions are written as 
$$
[A_{i+1}+A_i, B_{i+1}+B_i]=0\ \ {\rm and}\ \  [A_{i+1}-A_i, B_{i+1}-B_i]=0,
$$
or 
\begin{equation} \label{eq:ABi}
[A_i,B_i]+[A_{i+1},B_{i+1}]=0,\ [A_i,B_{i+1}]+[A_{i+1},B_i]=0.
\end{equation}

We will look for points $B_i$ written in two ways as follows:
$$
B_i=a_iA_{i-1}+b_iA_i=c_iA_i+d_iA_{i+1}.
$$
Then $B_{i+1}=a_{i+1}A_i+b_{i+1}A_{i+1}$. Taking cross-products with $A_i$ and $A_{i+1}$, we find that $c_i=b_{i+1}, d_i=a_{i+1}$. Since 
$$
[A_{i-1},A_i]=\ed_{2i-1}, [A_i,A_{i+1}]=\ed_{2i+1}, [A_{i-1},A_{i+1}]=\ve_{2i},
$$
we also have
$$
a_{i+1}\ed_{2i+1}=-a_i\ed_{2i-1},\ (b_{i+1}-b_i)\ed_{2i+1}=a_i\ve_{2i}.
$$

It follows that  $a_i = (-1)^i t/\ed_{2i-1}$
for some non-zero constant $t$ (non-zero since ${\bf B}$ is not homothetic to ${\bf A}$), and 
$$
b_{i+1}-b_i = t \frac{(-1)^i \ve_{2i}}{\ed_{2i-1}\ed_{2i+1}}.
$$
This system of linear equations on $b_i$ has a solution if and only if the sum of the right hand sides vanishes, as needed.
\end{proof}

\subsection{Bianchi permutability} \label{subsec:Bianchi}
 In the following theorem, the polygons are either closed or twisted.

\begin{theorem} \label{thm:Bianchi}
Assume that ${\bf P}\crel {\bf Q}$ and ${\bf P}\drel {\bf R}$. Then there exists a polygon ${\bf S}$ such that ${\bf Q}\drel {\bf S}$ and ${\bf R}\crel {\bf S}$.
\end{theorem}

\begin{proof}
Define point $S_1$ by requiring $P_1Q_1S_1R_1$ to be a centoraffine butterfly. 
Then $[R_1,S_1]=c, [Q_1,S_1]=d$. 
Let ${\mathcal L}_{Q_1 Q_2,d} (S_1)=S_2$ and ${\mathcal L}_{R_1 R_2,c} (S_1)=S'_2$. We claim that $S_2=S_2'$ and that $P_2Q_2S_2R_2$ is a centoaffine butterfly.

Indeed, one has
$$
{\mathcal L}_{S_1 R_1}(S'_2)=R_2,\ {\mathcal L}_{R_1 P_1}(R_2)=P_2,\ {\mathcal L}_{P_1 Q_1}(P_2)=Q_2,\ {\mathcal L}_{Q_1 S_1}(Q_2)=S_2.
$$
Since $P_1Q_1S_1R_1$ is a centoraffine butterfly, Lemma \ref{lm:butter} implies that the quadrilateral $S'_2R_2P_2Q_2 S_2$ closes up  and that it is a centoaffine butterfly. 

This shows that one can define the polygon ${\bf S}$ so that $P_iQ_iS_iR_i$ is a centoraffine butterfly for all $i$.
\end{proof}

\begin{figure}[ht]
\centering
\includegraphics[width=.8\textwidth]{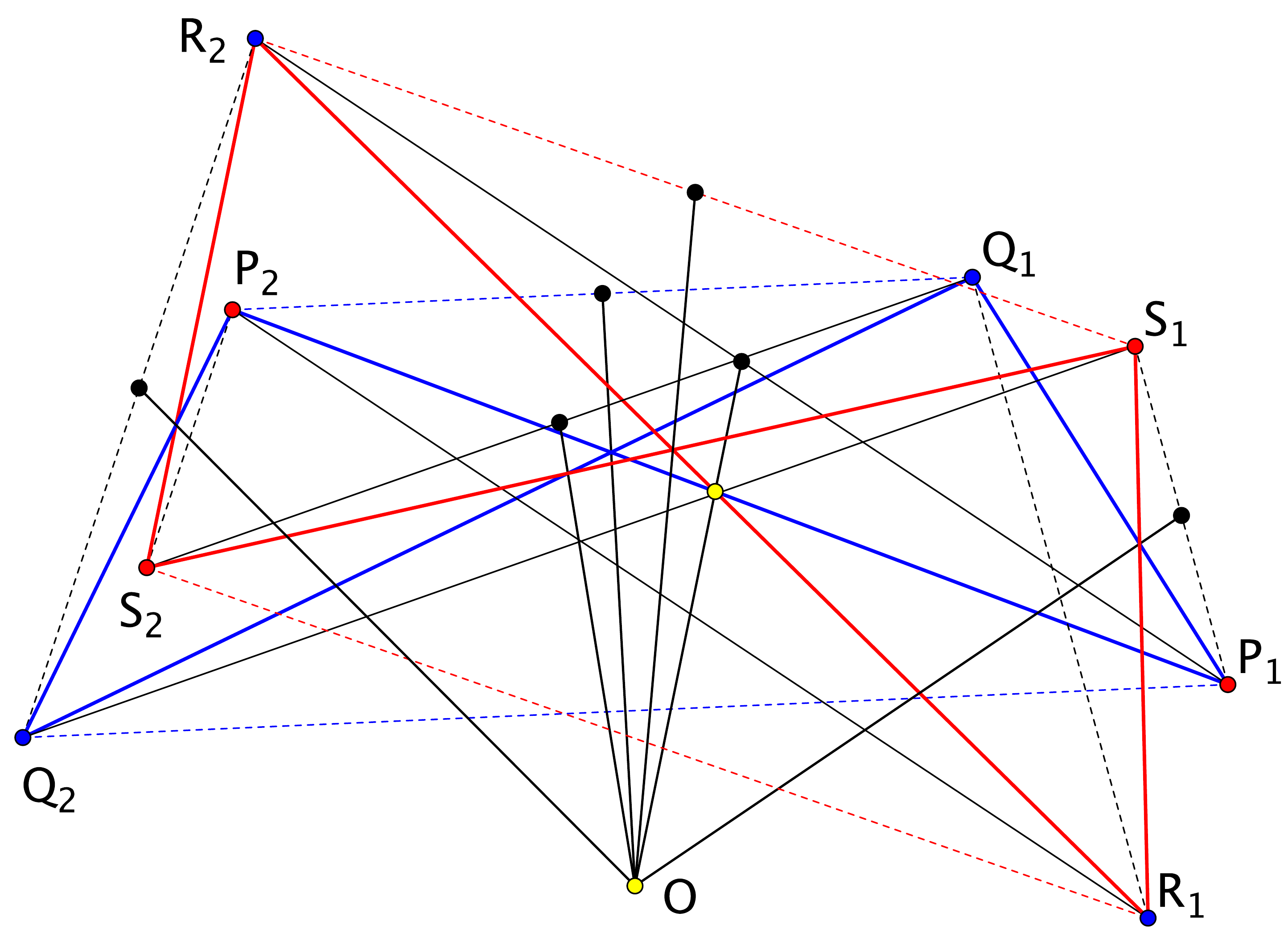}
\caption{Bianchi permutability.}
\label{Bianchi}
\end{figure}

We leave it to the reader to make sure that Theorem \ref{thm:Bianchi} can be interpreted as a configuration theorem depicted in Figure \ref{Bianchi}.

\section{Integrals} \label{sec:maps}

\subsection{Invariance of the Lax transformations} \label{subsec:inv}

The next result provide a Lax presentation of the $c$-relation. 

\begin{theorem} \label{thm:Lax}
If ${\bf P}$ and ${\bf Q}$ are $c$-related twisted $n$-gons, then the Lax transformations ${\mathcal L}_{{\bf P},\lambda}$ and ${\mathcal L}_{{\bf Q},\lambda}$ are conjugated for every value of the  parameter $\lambda$.
\end{theorem}

\begin{proof} 
According to Lemma \ref{lm:butter}, 
$$
{\mathcal L}_{Q_i Q_{i+1},\lambda} = {\mathcal L}_{P_{i+1} Q_{i+1},\lambda}\ {\mathcal L}_{P_i P_{i+1},\lambda}\ {\mathcal L}^{-1}_{P_{i} Q_{i},\lambda},
$$
and the monodromy acts trivially as a map $\R \to \R$.
Taking composition over $i=1,\ldots, n$ 
yields the result.
\end{proof}

Since ${\mathcal L}_{{\bf P},\lambda}$ is a fractional-linear transformation, one can realize it as a $2\times 2$ matrix, say, $A$. Then $({\rm Tr} A)^2/\det A$ is a well-defined conjugacy invariant function. It depends on $\lambda$, and expanding in a Taylor series in $\lambda$, one obtains integrals of the $c$-relation. 

\begin{example} \label{ex:tri}
{\rm  Consider the case of twisted triangles. Let $\ve_2,\ve_4,\ve_6$ and $\ed_{1},\ed_{3},\ed_{5}$ be the respective coordinates in the moduli space. 
Consider the space ${\mathcal Y}_{3,{\bf S}}$.

Multiply the three matrices (\ref{eq:Laxmat}) (replacing $c$ by $\lambda$) to find the Lax matrix ${\mathcal L}$ depending on this spectral parameter $\lambda$. One has
$$
{\Tr} {\mathcal L} = \lambda^3\left(\frac{\ve_2\ve_4\ve_6}{\ed_{1}\ed_{3}\ed_{5}}-\frac{\ve_2}{\ed_{5}}- \frac{\ve_4}{\ed_{1}} - \frac{\ve_6}{\ed_{3}} \right) +\lambda(\ve_2 \ed_{5} + \ve_4 \ed_{1} + \ve_6 \ed_{3}), 
$$
and
$$
\det {\mathcal L} = (\lambda^2-\ed_{1}^2) (\lambda^2-\ed_{3}^2) (\lambda^2-\ed_{5}^2).
$$
The determinant does not depend on the $\ve$-variables, hence the trace is invariant. We obtain two integrals:
$$
I=\frac{\ve_2\ve_4\ve_6}{\ed_{1}\ed_{3}\ed_{5}}-\left( \frac{\ve_2}{\ed_{5}}+ \frac{\ve_4}{\ed_{1}} + \frac{\ve_6}{\ed_{3}} \right)
$$
and
$$
J=\frac{\ve_2}{\ed_{1} \ed_3}+ \frac{\ve_4}{\ed_{3} \ed_5} + \frac{\ve_6}{\ed_{5} \ed_1}.
$$ 
Thus ${\mathcal Y}_{3,{\bf S}}$ is foliated by the common level curves of the functions $I$ and $J$.

For comparison, let ${\mathcal M}$ be the monodromy of the twisted triangle. Using formula (\ref{eq:monod}), we find that ${\Tr}\  {\mathcal M} =I$.
}
\end{example}

\begin{example} \label{ex:four}
{\rm 
Similarly, in the case of twisted quadrilaterals, we have variables $\ve_2,\ve_4,\ve_6,\ve_8$ and parameters $\ed_1,\ed_3,\ed_5,\ed_7$. The trace of the Lax matrix is a biquadratic polynomial in $\lambda$, and its three coefficients are integrals. The free term is $2\ed_1\ed_3\ed_5\ed_7$, a constant on ${\mathcal Y}_{4,{\bf S}}$, and the other two terms give two integrals
$$
I= \frac{\ve_2\ve_4\ve_6\ve_8}{\ed_1\ed_3\ed_5\ed_7} - 
\left( \frac{\ve_2\ve_4}{\ed_3\ed_7} + \frac{\ve_4\ve_6}{\ed_5\ed_1} + \frac{\ve_6\ve_8}{\ed_7\ed_3} + \frac{\ve_8\ve_2}{\ed_1\ed_5} \right) 
$$
and 
$$
J= \left( \frac{\ve_2\ve_4}{\ed_1\ed_3^2\ed_5} + \frac{\ve_4\ve_6}{\ed_3\ed_5^2\ed_7} + \frac{\ve_6\ve_8}{\ed_5\ed_7^2\ed_1} + \frac{\ve_8\ve_2}{\ed_7\ed_1^2\ed_3}\right) 
$$
(omitting the terms that do not depend on the $v$-variables).
Thus ${\mathcal Y}_{4,{\bf S}}$ is foliated by the common level surfaces of the functions $I$ and $J$. Note that again one has ${\Tr}\ {\mathcal M} = I$.

If the quadrilateral is closed, according to Example \ref{ex:34}, we have 
$$
\ve_6=-\ve_2, \ve_8=-\ve_4, \ve_2\ve_4=\ed_1\ed_5-\ed_3\ed_7,
$$
and both integrals, $I$ and $J$, become functions of the parameters $s_j$ only.
}
\end{example}

\begin{example} \label{ex:five}
{\rm Consider the case of twisted pentagons. As before, we calculate ${\Tr}\ {\mathcal L}$ and decompose it in homogeneous components in the spectral parameter $\lambda$. There are three terms, of degrees 1,3, and 5. This gives three integrals:
\begin{equation*}
\begin{split}
I=
\frac{\ve_2\ve_4\ve_6\ve_8\ve_{10}}{\ed_1\ed_3\ed_5\ed_7\ed_9} 
- &\left(\frac{\ve_2\ve_4\ve_6}{\ed_3\ed_5\ed_9} + \frac{\ve_4\ve_6\ve_8}{\ed_5\ed_7\ed_1} + \frac{\ve_6\ve_8\ve_{10}}{\ed_7\ed_9\ed_3} + \frac{\ve_8\ve_{10}\ve_{2}}{\ed_9\ed_1\ed_5} + \frac{\ve_{10}\ve_{2}\ve_{4}}{\ed_1\ed_3\ed_7}\right) \\
+ &\left( \frac{\ve_2\ed_7}{\ed_5\ed_9} + \frac{\ve_4\ed_9}{\ed_7\ed_1} + \frac{\ve_6\ed_1}{\ed_9\ed_3} + \frac{\ve_8\ed_3}{\ed_1\ed_5} + \frac{\ve_{10}\ed_5}{\ed_3\ed_7}  \right)
\end{split}
\end{equation*}
(this also comes from the trace of the monodromy),
\begin{equation*}
\begin{split}
&J=
\frac{\ve_2\ve_4\ve_6}{\ed_1\ed_3^2\ed_5^2\ed_7} + \frac{\ve_4\ve_6\ve_8}{\ed_3\ed_5^2\ed_7^2\ed_9} 
+ \frac{\ve_6\ve_8\ve_{10}}{\ed_5\ed_7^2\ed_9^2\ed_1} + \frac{\ve_8\ve_{10}\ve_{2}}{\ed_7\ed_9^2\ed_1^2\ed_3}
+ \frac{\ve_{10}\ve_{2}\ve_{4}}{\ed_9\ed_1^2\ed_3^2\ed_5}\\
&-  \frac{\ve_2}{\ed_1\ed_3} \left(\frac{1}{\ed_5^2}+\frac{1}{\ed_9^2}\right) - 
\frac{\ve_4}{\ed_3\ed_5} \left(\frac{1}{\ed_7^2}+\frac{1}{\ed_1^2}\right) -
\frac{\ve_6}{\ed_5\ed_7} \left(\frac{1}{\ed_9^2}+\frac{1}{\ed_3^2}\right) \\
&- \frac{\ve_8}{\ed_7\ed_9} \left(\frac{1}{\ed_1^2}+\frac{1}{\ed_5^2}\right) -
\frac{\ve_{10}}{\ed_9\ed_1} \left(\frac{1}{\ed_3^2}+\frac{1}{\ed_7^2}\right).
\end{split}
\end{equation*}
and
\begin{equation*}
\begin{split}
K=  \frac{\ve_2}{\ed_1\ed_3} + \frac{\ve_4}{\ed_3\ed_5} + \frac{\ve_6}{\ed_5\ed_7} + \frac{\ve_8}{\ed_7\ed_9} + \frac{\ve_{10}}{\ed_9\ed_1}.
\end{split}
\end{equation*}

The common level surfaces of these integrals foliate the 5-dimensional space ${\mathcal Y}_{5,{\bf S}}$.	
	
In the case of closed pentagons, we have the Ptolemy-Pl\"ucker relations
$$
\ve_2\ve_4=\ed_1\ed_5-\ed_3\ve_8
$$	
and its cyclic permutations. This makes the integrals functionally dependent, and leaves a single integral $K$ on the 2-dimensional  space ${\mathcal X}_{5,{\bf S}}$.
}
\end{example}

The formulas of examples \ref{ex:tri}, \ref{ex:four}, and \ref{ex:five} are extended to arbitrary values of $n$ in Sections \ref{subsec:inttwist} and \ref{subsec:intclosed} below.

\subsection{Centroaffine polygon recutting} \label{subsec:recut}

As mentioned in Introduction, by the {elementary recutting} of a closed polygon ${\bf P}$ we mean the linear transformation that changes only one vertex:
$P'_j=R_{P_{j-1} P_{j+1}} (P_j)$. Denote this transformation by $\rec_j$.	   Elementary recutings are involutions. 
The {recutting} $\rec$ of the polygon $\bf P$ is the composition $\rec_{n}\circ\cdots\circ\rec_{1}$.

The next lemma is a centroaffine analog of \cite[Cor. 2]{Ad2}.

\begin{lemma} \label{lm:Adler}
	Let $\rec_j$ be the elementary recutting at the vertex $j$. Then $(\rec_j \rec_{j+1})^3=Id$, $\rec_j^2=Id$ and $\rec_j\rec_k=\rec_k\rec_j$ for $|k-j|\geq 2$. 
\end{lemma}

\begin{figure}[!hbt]
	\centering
	\includegraphics[width=.9\textwidth]{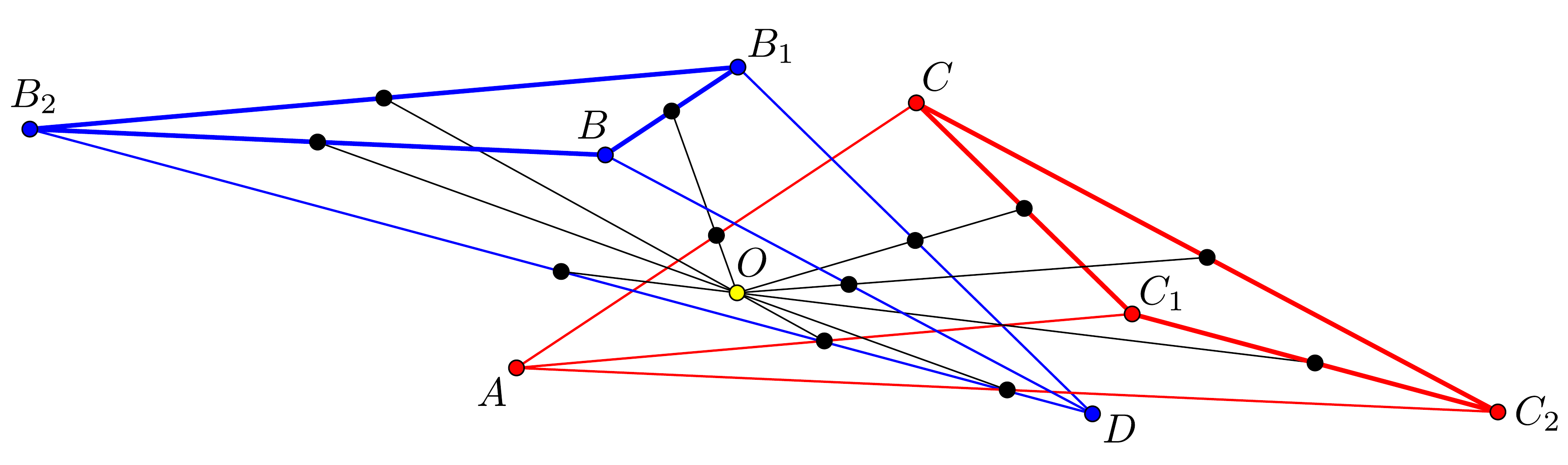}
	\caption{To Lemma \ref{lm:Adler}.}
	\label{Cube}
\end{figure}

\begin{proof}
The only non-trivial fact is the $(\rec_j \rec_{j+1})^3=Id$.
Consider Figure \ref{Cube} that depicts the configuration theorem described by this identity. One fixes the frame made of the vectors $A$ and $D$. Then we have six points 
$
B,B_1,B_2,C,C_1,C_2
$	
that should satisfy 12 relations (cf. (\ref{eq:cdef})):
\begin{equation*}
\begin{split}
&[A,B]=[B_1,C], [A,B_1]=[B,C], [A,B_1]=[B_2,C_1], [A,B_2]=[B_1,C_1],\\
&[A,B_2]=[B,C_2], [A,B]=[B_2,C_2], [C_2,D]=[B,C], [C,D]=[B,C_2],\\
&[C,D]=[B_1,C_1], [C_1,D]=[B_1,C], [C_1,D]=[B_2,C_2], [C_2,D]=[B_2,C_1].
\end{split}
\end{equation*}
This system is equivalent to the system of 9 relations
\begin{equation} \label{eq:nine}
\begin{split}
&[A,B]=[B_1,C]=[B_2,C_2]=[C_1,D] :=x,\\
&[A,B_1]=[B_2,C_1]=[B,C]=[C_2,D] :=y,\\
&[A,B_2]=[B,C_2]=[B_1,C_1]=[C,D] :=z.		
\end{split}
\end{equation}	
Assume, without loss of generality, that $[A,D]=1$ and express the remaining points as linear combinations of $A$ and $D$:
\begin{equation*}
\begin{split}
B=\lambda A+xD, B_1=\lambda_1 A + yD, B_2=\lambda_2 A + zD,\\
C=\mu D +zA, C_1=\mu_1 D + xA, C_2=\mu_2 D + yA.
\end{split}
\end{equation*}
Then the rest of the equations (\ref{eq:nine}) become the following six equations on the nine variables $x,y,z,\lambda_i, \mu_i, i=1,2,3$:
$$
\lambda_1\mu=x+yz=\lambda_2\mu_2, \lambda_2\mu_1=y+xz=\lambda\mu, \lambda\mu_2=z+xy=\lambda_1\mu_1.
$$
These equations are not independent: the product of the three left hand sides equals the product of the three right hand sides. 

Overall, one has nine variables satisfying five relations. The resulting four degrees of freedom make it possible to choose point $B$ and $C$ arbitrarily, proving the existence of the configuration of the lemma.
\end{proof}

\begin{remark}
{\rm The relations of Lemma \ref{lm:Adler} show that one has a representation of the  group of permutations $S_{n+1}$ on the space of centroaffine $n$-gons, similarly to polygon recutting, see \cite{Ad2}.
}
\end{remark}

The next lemma states that elementary recuttings commute with the $c$-transformations (an analog of one of the statements of Theorem 4 in \cite{TT}).

\begin{lemma}\label{lm:elementary_recutting}
If ${\bf P}\crel {\bf Q}$ then $\rec_j({\bf P})\crel \rec_j({\bf Q})$ for all $j$.
\end{lemma}

\begin{proof}
	One only has to check that $[P'_j,Q'_j]=c$. 
	From \eqref{eq:P_j+1} it follows that
	$$P_j=\frac{\ed_{2j+1}({\bf{P}})}{\ve_{2j}({\bf{P}})}P_{j-1}+\frac{\ed_{2j-1}({\bf{P}})}{\ve_{2j}({\bf{P}})} P_{j+1}.
	$$
	Since $[P_j',P_{j+1}]=[P_{j-1},P_j]$ and $[P_{j-1},P_j']=[P_j,P_{j+1}]$, one has \begin{equation}\label{eq:elem_recut}P_j'=\frac{\ed_{2j-1}({\bf{P}})}{\ve_{2j}({\bf{P}})} P_{j-1}+\frac{\ed_{2j+1}({\bf{P}})}{\ve_{2j}({\bf{P}})} P_{j+1}.\end{equation}
	
Likewise, for $\bf Q$ one has
$$
Q_j=\frac{\ed_{2j+1}({\bf{Q}})}{\ve_{2j}({\bf{Q}})}Q_{j-1}+\frac{\ed_{2j-1}({\bf{Q}})}{\ve_{2j}({\bf{Q}})} Q_{j+1},
$$
and
$$
Q_j'=\frac{\ed_{2j-1}({\bf{Q}})}{\ve_{2j}({\bf{Q}})} Q_{j-1}+\frac{\ed_{2j+1}({\bf{Q}})}{\ve_{2j}({\bf{Q}})}Q_{j+1}.
$$
	Since ${\bf Q}\crel {\bf P}$ we have
\begin{equation*}
\begin{split}
c&=[P_j,Q_j]=
\left(\tfrac{\ed_{2j+1}({\bf{P}})}{\ve_{2j}({\bf{P}})}\tfrac{\ed_{2j+1}({\bf{Q}})}{\ve_{2j}({\bf{Q}})}+\tfrac{\ed_{2j-1}({\bf{P}})}{\ve_{2j}({\bf{P}})}\tfrac{\ed_{2j-1}({\bf{Q}})}{\ve_{2j}({\bf{Q}})}\right)c\\
&+\tfrac{\ed_{2j+1}({\bf{P}})}{\ve_{2j}({\bf{P}})}\tfrac{\ed_{2j-1}({\bf{Q}})}{\ve_{2j}({\bf{Q}})}[P_{j-1},Q_{j+1}]+\tfrac{\ed_{2j-1}({\bf{P}})}{\ve_{2j}({\bf{P}})}\tfrac{\ed_{2j+1}({\bf{Q}})}{\ve_{2j}({\bf{Q}})}[P_{j+1},Q_{j-1}].
\end{split}
\end{equation*}
	For $[P_j',Q_j']$, we have
\begin{equation*}
\begin{split}
[P_j',Q_j']&=\left(\tfrac{\ed_{2j+1}({\bf{P}})}{\ve_{2j}({\bf{P}})}\tfrac{\ed_{2j+1}({\bf{Q}})}{\ve_{2j}({\bf{Q}})}+\tfrac{\ed_{2j-1}({\bf{P}})}{\ve_{2j}({\bf{P}})}\tfrac{\ed_{2j-1}({\bf{Q}})}{\ve_{2j}({\bf{Q}})}\right)c\\
&+\tfrac{\ed_{2j-1}({\bf{P}})}{\ve_{2j}({\bf{P}})}\tfrac{\ed_{2j+1}({\bf{Q}})}{\ve_{2j}({\bf{Q}})}[P_{j-1},Q_{j+1}]+\tfrac{\ed_{2j+1}({\bf{P}})}{\ve_{2j}({\bf{P}})}\tfrac{\ed_{2j-1}({\bf{Q}})}{\ve_{2j}({\bf{Q}})}[P_{j+1},Q_{j-1}].
\end{split}
\end{equation*}
	But 
\begin{equation*}
\begin{split}
	\frac{\ed_{2j-1}({\bf{P}})}{\ve_{2j}({\bf{P}})}\frac{\ed_{2j+1}({\bf{Q}})}{\ve_{2j}({\bf{Q}})}&=\frac{[P_{j-1},P_{j}][Q_{j},Q_{j+1}]}{[P_{j-1},P_{j+1}][Q_{j-1},Q_{j+1}]}\\
	&=\frac{[Q_{j-1},Q_{j}][P_{j},P_{j+1}]}{[P_{j-1},P_{j+1}][Q_{j-1},Q_{j+1}]}=\frac{\ed_{2j+1}({\bf{P}})}{\ve_{2j}({\bf{P}})}\frac{\ed_{2j-1}({\bf{Q}})}{\ve_{2j}({\bf{Q}})}.
\end{split}
\end{equation*}
	 Hence the result.
	 \end{proof}

\begin{figure}[ht]
\centering
\includegraphics[width=.65\textwidth]{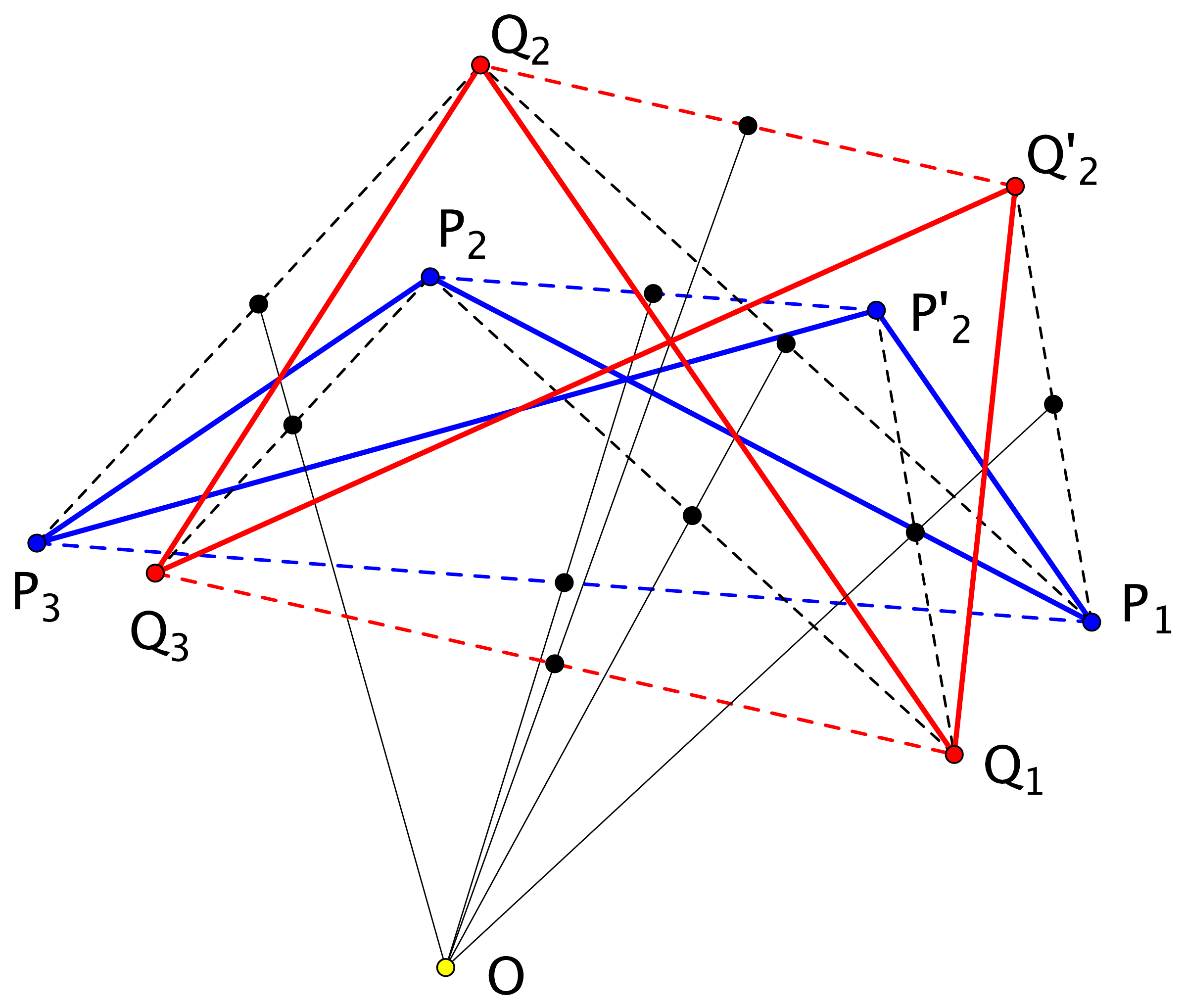}
\caption{Elementary recutting commutes with $c$-relation.}
\label{recutting}
\end{figure}

Lemma \ref{lm:elementary_recutting}  also can be interpreted as a configuration theorem, see Figure \ref{recutting}. As before, we leave it to the reader to make sure of it.

\begin{theorem} \label{thm:recut}
1) The Lax transformation is preserved by recutting;\\
2) Recutting commutes with $c$-relations.	
\end{theorem}

\begin{proof}
The first statement follows from Lemma \ref{lm:butter}:
$$
{\mathcal L}_{P_i P_{i+1}} {\mathcal L}_{P_{i-1} P_{i}}  = {\mathcal L}_{P'_i P_{i+1}} {\mathcal L}_{P_{i-1} P'_{i}} 
$$
for all $i$.
The second statement follows from Lemma \ref{lm:elementary_recutting}: if $\bf P\crel Q$ then $\rec_1({\bf P})\crel \rec_1({\bf Q})$, therefore $\rec_2\rec_1({\bf P})\crel \rec_2\rec_1({\bf Q})$, etc. 
\end{proof}

\subsection{Odd-gons, infinitesimal map} \label{subsec:infin}

If the constant $c$ is infinitesimal, we obtain an $SL(2,\R)$-invariant vector field on the space of twisted polygons, that is, a vector field on the moduli space ${\mathcal Y}_n$. In this section we calculate this field. 

Let ${\mathbf P}$ be a twisted $n$-gon with odd $n$.
Let the field be given by vectors $\xi_i$ with foot points $P_i,\ i=1,\ldots,n$. The conditions $[P_i,Q_i]=c, [P_i,P_{i+1}]=[Q_i,Q_{i+1}]$ become
\begin{equation} \label{eq:field1}
[P_i,\xi_i]=1, \ [P_i,\xi_{i+1}]+[\xi_i,P_{i+1}]=0
\end{equation}
(we normalize the field so the constant in the first equation is 1).

\begin{theorem} \label{thm:field}
In terms of the $(v,s)$-coordinates, the field $\xi$ is given by
$$
\dot v_{2i} = v_{2i} \left(\sum_{k=1}^{n-1} (-1)^{k-1} \frac{v_{2i+2k}}{\ed_{2i+2k-1}\ed_{2i+2k+1}}\right) +\frac{\ed_{2i+1}}{\ed_{2i-1}}-\frac{\ed_{2i-1}}{\ed_{2i+1}}.
$$
\end{theorem}


\begin{proof}
Let $\xi_i=a_iP_i+b_i P_{i+1}$. Then (\ref{eq:field1}) implies
$$
b_i \ed_{2i+1}=1, (a_i+a_{i+1})\ed_{2i+1} + b_{i+1} v_{2i+2} =0.
$$
Since $n$ is odd, there is a unique solution
\begin{equation} \label{eq:field2}
a_i=\frac{1}{2} \sum_{k=1}^{n} (-1)^k \frac{v_{2i+2k}}{\ed_{2i+2k-1}\ed_{2i+2k+1}}, \ b_i=\frac{1}{\ed_{2i+1}}.
\end{equation}

Since $v_{2i}=[P_{i-1},P_{i+1}]$, we have 
\begin{equation} \label{eq:field3}
\begin{split}
&\dot v_{2i}=[P_{i-1},\xi_{i+1}]+[\xi_{i-1},P_{i+1}] 
= [P_{i-1},a_{i+1}P_{i+1}+b_{i+1} P_{i+2}]\\
&+[a_{i-1}P_{i-1}+b_{i-1} P_{i}, P_{i+1}]
= (a_{i-1}+a_{i+1})v_{2i} + b_{i+1} [P_{i-1},P_{i+2}]+ b_{i-1}\ed_{2i+1}.
\end{split}
\end{equation} 
By the Prolemy-Pl\"ucker relation, 
$$
[P_{i-1},P_{i+2}]=\frac{v_{2i}v_{2i+2}-\ed_{2i-1}\ed_{2i+3}}{\ed_{2i+1}},
$$
and the values of $a_i$ and $b_i$ are found in (\ref{eq:field2}).  Substitute in (\ref{eq:field3})
to obtain the result.
\end{proof}

It follows from Theorems \ref{thm:Bianchi} and \ref{thm:recut} that the flow of the field $\xi$ commutes with the $c$-relations and with the polygon recutting.

\begin{remark}
{\rm If $n$ is even, then a necessary condition for $a_i$ to exist is 
$$
\sum_{i=1}^n (-1)^n \frac{v_{2i}}{\ed_{2i-1}\ed_{2i+1}}=0,
$$
the equation that appeared in Lemma \ref{lm:ident}.
If this necessary condition holds, the vector field is defined modulo a 1-dimensional space, the kernel of the matrix that gives the linear equations on the variables $a_i$.
}
\end{remark}

Set
$$
g_i = \frac{v_{2i}}{\ed_{2i-1}\ed_{2i+1}},\ \beta_i=-\frac{1}{\ed_{2i-1}^2}.
$$
Then, after scaling, the vector field $\xi$ is given by
$$
\dot g_i=-g_i(g_{i+1}-g_{i+2}+g_{i+3}-\ldots -g_{i+n-1}) +\beta_i-\beta_{i+1}.
$$
This is the dressing chain of Veselov-Shabat, formula (12) in \cite{VS}. 

\subsection{Integrals for twisted polygons} \label{subsec:inttwist}

In view of Theorem \ref{thm:Lax}, we calculate the trace of the Lax matrix ${\mathcal L}$ of an $n$-gon, the product of $n$ matrices (\ref{eq:Laxmat}). Notice that 
$$
\begin{pmatrix}
	0&1\\
	1&0
\end{pmatrix}
\begin{pmatrix}
	0&1\\
	-b&a
\end{pmatrix}
\begin{pmatrix}
	0&1\\
	1&0
\end{pmatrix}^{-1}
=
\begin{pmatrix}
	a&-b\\
	1&0
\end{pmatrix}.
$$
Therefore we can use the calculations from Section \ref{subsec:mono}.  

Set 
\begin{equation} \label{eq:abi}
a_i=\frac{\lambda\ve_{2i}}{\ed_{2i-1}\ed_{2i+1}},\ b_i=\frac{\lambda^2}{\ed^2_{2i-1}}-1,\ i=0,\ldots,n-1,
\end{equation}
and consider the continuants (\ref{eq:continuant_ij}). Then Proposition \ref{prop:monod} implies that
$$
{\rm Tr}\ {\mathcal L} = D_{0,n}-b_0D_{1,n-1}.
$$
The homogeneous in $\lambda$ components are integrals of the $c$-relations for all values of $c$.

A combinatorial rule for calculation of the general continuants 
\[
\begin{vmatrix}
a_1 & c_1 & 0 & \ldots & 0 & 0\\
b_1 & a_2 & c_2 & \ldots & 0 & 0\\
0 & b_2 & a_3 & \ldots & 0 & 0\\
\vdots & \vdots & \vdots & \ddots & \vdots & \vdots\\
0 & 0 & 0 & \ldots &a_{n-1} & c_{n-1}\\
0 & 0 & 0 & \ldots &b_{n-1} & a_n
\end{vmatrix}
\]
is as follows: one term of the continuant is $a_1a_2\ldots a_n$, and the other terms are obtained from it by replacing any number of disjoint pairs $(a_i a_{i+1})$ by $(-b_ic_i)$, see \cite{Muir}. That is, this continuant can be written as
\begin{equation} \label{eq:gendiff}
\prod_{i=1}^{n-1} \left(E-b_ic_i\frac{\partial^2}{\partial a_i \partial a_{i+1}}\right) (a_1a_2\cdots a_n),
\end{equation}
where $E$ is the identity operator.
Note that the differential operators involved in this formula commute with each other.

As in \cite{AFIT}, a subset of the set $\{0,1,\ldots,n-1\}$ is called {\it cyclically sparse} if it contains no pairs of consecutive indices, and the indices are understood cyclically mod $n$ (the empty set is also sparse). 

The above rule implies 

\begin{lemma} \label{lm:genfunc}
One has
\begin{equation} \label{eq:genfunc}
\begin{split}
{\rm Tr}\ {\mathcal L} &= \sum_{I\ {\rm sparse}} (-1)^{|I|}  \lambda^{n-2|I|}   \prod_{j,j+1\notin I} \left(\frac{\ve_{2j}}{\ed_{2j-1}\ed_{2j+1}}\right) \prod_{i \in I} \left( \frac{\lambda^2}{\ed^2_{2i-1}}-1 \right)\\
&= \lambda^n \sum_{I\ {\rm sparse}} \prod_{j,j+1\notin I} \left(\frac{\ve_{2j}}{\ed_{2j-1}\ed_{2j+1}}\right) \prod_{i \in I} \left(\frac{1}{\lambda^{2}}-\frac{1}{\ed^2_{2i-1}}\right).
\end{split}
\end{equation}
\end{lemma}

Thus, for  $t=\lambda^{-2}$, 
\begin{equation} \label{eq:integr}
\sum_{I\ {\rm sparse}} \prod_{j,j+1\notin I} \left(\frac{\ve_{2j}}{\ed_{2j-1}\ed_{2j+1}}\right) \prod_{i \in I} \left(\frac{1}{\lambda^{2}}-\frac{1}{\ed^2_{2i-1}}\right)=: \sum_k t^k F_k,
\end{equation}
is a generating function of the integrals $F_k$ of the $c$-relation on twisted $n$-gons.

Note that, as a polynomial in the variables $\ve_{2i}$, one has deg $F_k=n-2k$. In particular, if $n=2q+1$, then
$$
F_q=\sum_i \frac{\ve_{2i}}{\ed_{2i-1}\ed_{2i+1}},
$$
and one has $q+1$ integrals in this case. If $n=2q$, then 
$$
F_{q-1}= \sum_i \frac{\ve_{2i}\ve_{2i+2}}{\ed_{2i-1}\ed_{2i+1}^2\ed_{2i+3}},
$$
and one has $q$ integrals in this case. Compare with the examples in Section \ref{subsec:inv}.

To summarize, we obtain $\lfloor{\frac{n+1}{2}}\rfloor$ integrals on the moduli space of twisted $n$-gons. We do not dwell on the independence and completeness of this set of integrals here; see \cite{Iz} and the remark below.

\begin{remark} 
{\rm If one has $s_{2i-1}=1$ for all $i$, corresponding -- as explained in Introduction -- to the case  studied in \cite{AFIT}, then formula (\ref{eq:genfunc})  simplifies. The integrals in \cite{AFIT} are given by the formulas
$$
F_k = \sum_{I\ {\rm sparse}, |I|=k} \prod_{i\in I} c_i,
$$
where, in our notation, $c_i=1/(\ve_{2i}\ve_{2i+2})$. Multiplying by $(\prod c_i)^{-1/2}$ (which is also an integral), this becomes 
$$
\sum_{I\ {\rm sparse}, |I|=k} \prod_{j,j+1\notin I} \ve_{2j},
$$
precisely what (\ref{eq:genfunc}) yields when $s_{2i-1}=1$ for all $i$.
}
\end{remark}

As before, we write
$$
g_i = \frac{v_{2i}}{\ed_{2i-1}\ed_{2i+1}},\ \beta_i=-\frac{1}{\ed_{2i-1}^2}.
$$
Modifying formula (\ref{eq:gendiff}) to account for the cyclic symmetry, we 
write the generating function of the integrals $F_k$ as
$$
 \prod_{i=0}^{n-1} \left[ E+ (t+\beta_{i+1}) \frac{\partial^2}{\partial g_i \partial g_{i+1}} \right]
\left(\prod_{j=0}^{n-1} g_j \right).
$$
Comparing this formula with formula (22) in \cite{VS}, we conclude that our integrals coincide with the integrals of the dressing chain  obtained \cite{VS}. This is not unexpected, given the appearance of the dressing chain as the infinitesimal version of the $c$-relation in the preceding  Section \ref{subsec:infin}. Note that the independence of the integrals is asserted in \cite{VS}.

\subsection{Integrals for closed polygons} \label{subsec:intclosed}

When restricted to the moduli space of closed polygons, the integrals from Section \ref{subsec:inttwist} become dependent. This follows from the next general observation.

Let $X$ be a manifold, $M_t: X \to SL(2,\R)$ be  a 1-parameter family of smooth maps, analytically depending on parameter $t$. Assume that the unit matrix $E\in SL(2,\R)$ is a regular value of $M_0$, and let $Y=M_0^{-1} (E)$. Consider a 1-parameter family of smooth functions  $G_t={\Tr}M_t: X \to \R$. Let prime denote $d/dt$.

\begin{lemma} \label{lm:Tayl}
One has: 
$$
G_0\left|_Y\right.=2,\ G'_0\left|_Y\right.=0,\ dG_0\left|_Y\right.=0.
$$
\end{lemma}

\begin{proof}
Since $M_0$ sends $Y$ to the unit matrix, $G_0$ equals 2 on $Y$. 

Fix $y\in Y$ and consider the curve $M_t(y)$ in $SL(2,\R)$. The tangent vector $M'_t (y)$ at $t=0$ lies in $sl(2,\R)$, and this matrix has zero trace. This proves the second equality. 

For the third equality, let the eigenvalues of $M_t(x)$ be $e^{\pm \mu(x,t)}$. Then $G_t(x)=e^{\mu(x,t)} + e^{-\mu(x,t)}$, hence
$dG_t = (e^{\mu(x,t)} - e^{-\mu(x,t)}) d\mu(x,t)$. Since $\mu(x,0)=0$ for $x\in Y$, one has $dG_0\left|_Y\right.=0$, as claimed.
\end{proof}

We shall apply this lemma to a modified Lax matrix of an $n$-gon. Recall that $t=\lambda^{-2}$. Let
\begin{equation} \label{eq:barL}
\bar{\mathcal L}=\frac{{\mathcal L}}{\sqrt{\prod \left(\frac{\lambda^2}{\ed_{2i-1}^2} -1 \right)}} =  \frac{ {\mathcal L} \left(\prod_{i=1}^n s_{2i-1}\right)}{\lambda^n\sqrt{\prod \left(1- t\ed_{2i-1}^2 \right)}}. 
\end{equation}
Then, according to (\ref{eq:Laxmat}), $\bar{\mathcal L} \in SL(2,\R)$. We write $\bar{\mathcal L}_t$ to emphasize the dependence on $t$.

\begin{lemma}  \label{lm:Laxclosed}
One has $\bar{\mathcal L}_0=E$ on closed polygons. 
\end{lemma}

\begin{proof}
We will show that the Lax matrix $\bar{\mathcal L}_t$ is a deformation of the monodromy ${\mathcal M}$. 

The monodromy was calculated in Proposition (\ref{prop:monod}),  where the variables in the continuants were
$$
a_i=\frac{v_{2i}}{s_{2ji1}},\ b_i=\frac{s_{2i+1}}{s_{2i-1}}.
$$
One can rewrite these continuants as follows: 
$$
D_{i,j+1}=	(s_{2i+1} s_{2i+3}\cdots s_{2j+1}) \bar D_{i,j+1},
$$
 where
$$
\bar D_{i,j+1}=	
\det     
\begin{pmatrix}
	\frac{\ve_{2i}}{\ed_{2i-1}\ed_{2i+1}}&1&0&\cdots&0\\
	\frac{1}{\ed^2_{2i+1}}&\frac{\ve_{2i+2}}{\ed_{2i+1}\ed_{2i+3}}&1&0&\cdots\\
	0&\ddots&\ddots&\ddots&0\\
	\cdots&0&\frac{1}{\ed^2_{2j-3}}&\frac{\ve_{2j-2}}{\ed_{2j-3}\ed_{2j-1}}&1\\
	0&\cdots&0&\frac{1}{\ed^2_{2j-1}}&\frac{\ve_{2j}}{\ed_{2j-1}\ed_{2j+1}}
\end{pmatrix}.
$$
Since ${\mathcal M}=E$ for closed $n$-gons, we have
\begin{equation} \label{eq:dbar}
-\left(\frac{\prod_{i=1}^n \ed_{2i-1}}{\ed^2_{2n-1}}\right) \bar D_{1,n-1} = \left(\prod_{i=1}^n \ed_{2i-1}\right) \bar D_{0,n}=1,\ \bar D_{1,n}=\bar D_{0,n-1}=0.
\end{equation}

On the other hand, ${\mathcal L}_0$ is also given by Proposition (\ref{prop:monod}),  where the variables in the continuants are
as in (\ref{eq:abi}). This matrix equals
$$
\lambda^n 
\begin{pmatrix}
-\left(\frac{1}{\ed^2_{2n-1}}-t \right) \bar D_{1,n-1}(t),& \frac{1}{\lambda} \bar D_{0,n-1}(t)\\
-\lambda \left(\frac{1}{\ed^2_{2n-1}}-t \right) \bar D_{1,n}(t),&\bar D_{0,n}(t)
\end{pmatrix},
$$
where 
$$
\bar D_{i,j+1}(t)=	
\det     
\begin{pmatrix}
	\frac{\ve_{2i}}{\ed_{2i-1}\ed_{2i+1}}&1&0&\cdots&0\\
	\frac{1}{\ed^2_{2i+1}}-t&\frac{\ve_{2i+2}}{\ed_{2i+1}\ed_{2i+3}}&1&0&\cdots\\
	0&\ddots&\ddots&\ddots&0\\
	\cdots&0&\frac{1}{\ed^2_{2j-3}}-t&\frac{\ve_{2j-2}}{\ed_{2j-3}\ed_{2j-1}}&1\\
	0&\cdots&0&\frac{1}{\ed^2_{2j-1}}-t&\frac{\ve_{2j}}{\ed_{2j-1}\ed_{2j+1}}
\end{pmatrix}.
$$

Since $\bar D_{1,n}=0$, one has $\bar D_{1,n}(t)=O(t)$, and hence  $\lim_{t\to 0} \lambda \bar D_{1,n}(t) =0.$ This, along with equations (\ref{eq:dbar}), implies that $\bar{\mathcal L}_0=E$, as claimed.
\end{proof}

We apply Lemma \ref{lm:Tayl} to $\bar{\mathcal L}_t$. Recall that ${\Tr} {\mathcal L}_t = \lambda^n (F_0+tF_1+\ldots)$ where $F_i$ are  integrals of the $c$-relation on twisted $n$-gons.

\begin{proposition} \label{prop:relint}
Restricted to the space ${\mathcal X}_n$ of closed $n$-gons, one has the following relations:
$$
F_0\left|_{\mathcal X}\right. =2\left(\prod_{i=1}^n \ed_{2i-1}\right)^{-1/2},\ F_1\left|_{\mathcal X}\right. = -\left(\frac{1}{2} \sum_{i=1}^n \ed_{2i-1}\right) F_0\left|_{\mathcal X}\right.,\ dF_0\left|_{\mathcal X}\right.=0.
$$
\end{proposition}

\begin{proof}
One has 
$$
\prod_{i=1}^n \left(1- t\ed_{2i-1}\right)^{-1/2} = 1 + \left(\frac{1}{2} \sum_{i=1}^n \ed_{2i-1}\right) t + O(t^2),
$$
and then, according to (\ref{eq:barL}), 
\begin{equation*}
\begin{split}
{\Tr}\bar{\mathcal L}_t &= \left(\prod_{i=1}^n \ed_{2i-1}\right) \left[1+\left(\frac{1}{2} \sum_{i=1}^n \ed_{2i-1}\right) t + O(t^2)\right] \left[F_0+tF_1+O(t^2)  \right]\\
&= \left(\prod_{i=1}^n \ed_{2i-1}\right) \left[ F_0 + \left(F_1 + \left(\frac{1}{2} \sum_{i=1}^n \ed_{2i-1}\right) F_0 \right) t + O(t^2)     \right].
\end{split}
\end{equation*}
Now Lemma \ref{lm:Tayl} implies the result.
\end{proof}

As before, we do not dwell here on the question whether the relations from Proposition \ref{prop:relint} are the only ones satisfied by the integrals when restricted to the moduli space of closed polygons (similarly to \cite{AFIT}, we do expect this to be the case).

\section{Closed centroaffine  polygons, before factorizing by $SL(2,\R)$} \label{sec:before}

\subsection{Presymplectic forms} \label{subsec:forms}
Recall that $\widetilde{\mathcal X}_{n,{\bf S}}$ is the space of closed $n$-gons with fixed  ``side areas" $s_{2i-1},\ i=1,\ldots,n$.

Choose a coordinate system in $\R^2$, and let $P_i=(x_i,y_i)$ be the vertices of an $n$-gon. Consider the differential 2-form
\begin{equation} \label{eq:omega}
\omega = \sum_{i=1}^n s_{2i+1} (dx_i \wedge dy_{i+1} + dx_{i+1} \wedge dy_{i})
\end{equation}
in $\widetilde{\mathcal X}_{n}$. The restriction of $\omega$ to $\widetilde{\mathcal X}_{n,{\bf S}}$, which we denote by the same letter, is closed, that is,  a presymplectic form in $\widetilde{\mathcal X}_{n,{\bf S}}$. 

The generators of the Lie algebra $sl(2,\R)$ acting diagonally on polygons are the vector fields
$$
e=\sum_i x_i \partial/\partial y_i,\ h=\sum_i (x_i \partial/\partial x_i - y_i \partial/\partial y_i),\ f=\sum_i y_i \partial/\partial x_i. 
$$
These vector fields are tangent to the submanifolds $\widetilde{\mathcal X}_{n,{\bf S}}$.

Let
$$
I= \sum_i s_{2i+1} x_i x_{i+1},\  J= \sum_i s_{2i+1} (x_i y_{i+1} + x_{i+1}y_i),\ K= \sum_i s_{2i+1} y_iy_{i+1}. 
$$
The restrictions of these functions on $\widetilde{\mathcal X}_{n,{\bf S}}$ are the Hamiltonian functions of the above vector fields:
\begin{equation} \label{eq:Ham}
i_e \omega=-dI,\ i_h \omega=dJ,\ i_f\omega=dK.
\end{equation}

\begin{theorem} \label{thm:pres}
1) The restriction of $\omega$ to $\widetilde{\mathcal X}_{n,{\bf S}}$
is $SL(2,\R)$-invariant, but it is not basic: it does not descend on the moduli space ${\mathcal X}_{n,{\bf S}}$.\\
2) The form $\omega$ is invariant under the $c$-relations for all values of $c$  and under the polygon recutting.
\end{theorem}


\begin{proof}
Equations (\ref{eq:Ham}) and the Cartan formula imply that $L_e(\omega)=L_h(\omega)=L_f(\omega)=0$. Therefore $\omega$ is $SL(2,\R)$-invariant, but it is not basic since it is not annihilated by $sl(2,\R)$.

To prove the invariance of $\omega$ under the $c$-relations, let ${\bf Q}$ be an  $n$-gon such that $\mathbf{Q}\crel \mathbf{P}$, and let $Q_i=(u_i,v_i)$. One has
$$
s_{2i+1} Q_i = c P_{i+1} + [Q_i,P_{i+1}] P_i,\
s_{2i+1} P_{i+1} = c Q_i + [Q_i,P_{i+1}] Q_{i+1}.
$$
Take bracket of the first equality with $dQ_{i+1}$, bracket of the second equality with $dP_i$, subtract the second from the first, and sum up over $i$:
\begin{equation} \label{eq:1form}
\begin{split}
\sum s_{2i+1} ([Q_i,dQ_{i+1}] - [P_{i+1},dP_i]) = c \sum ([P_{i+1},dQ_{i+1}] - [Q_i,dP_i]) \\
+ \sum ([Q_i,P_{i+1}] ([P_i,dQ_{i+1}] - [Q_{i+1},dP_i])).
\end{split}
\end{equation}
The differential of the left hand side of (\ref{eq:1form}) is the difference of $\omega$ evaluated at ${\bf Q}$ and ${\bf P}$. Therefore it suffices to show that the right hand side is an exact 1-form on $\widetilde{\mathcal X}_{n,{\bf S}}$.

Indeed, 
$$
\sum ([P_{i+1},dQ_{i+1}] - [Q_i,dP_i]) = \sum ([P_{i},dQ_{i}] - [Q_i,dP_i]) = \sum d(x_iv_i+y_iu_i).
$$
Next, $[P_i,dQ_{i+1}] - [Q_{i+1},dP_i] = d [P_i,Q_{i+1}]$. Equation (\ref{eq:invol}) implies
$$
[P_i,Q_{i+1}] = \frac{s_{2i+1}^2-c^2}{[Q_i,P_{i+1}]}.
$$
It follows that, up to a constant, 
$$
[Q_i,P_{i+1}] ([P_i,dQ_{i+1}] - [Q_{i+1},dP_i]) =  [Q_i,P_{i+1}] d([Q_i,P_{i+1}]^{-1}) = - d \ln([Q_i,P_{i+1}]), 
$$
therefore the second sum on the right hand side of (\ref{eq:1form}) is also an exact 1-form.

To show that the form $\omega$ is invariant under the polygon recutting, consider the difference of the form evaluated at  $\mathbf{P}$ and at $\rec_{1}(\mathbf{P})$. Let $P'_1=(x,y)=\rec_{1}({\bf P})_1$. Then 
$$
x=\dfrac{x_n \ed_1+x_2 \ed_3}{\ve_2},\ y=\dfrac{y_n \ed_1+y_2 \ed_3}{\ve_2}. 
$$
We get
\begin{equation*}
	\begin{split}
		\omega(\mathbf{P})-\omega(\rec_1 (\mathbf{P}))&=\ed_3(dx_1\wedge dy_2+dx_2\wedge dy_1-dx_n\wedge dy-dx\wedge dy_n)\\
		&+\ed_1(dx_n\wedge dy_1+dx_1\wedge dy_n-dx\wedge dy_2-dx_2\wedge dy)\\
		=(\ed_3dx_1-\ed_1 dx)&\wedge dy_2+dx_2\wedge(\ed_3 dy_1-\ed_1 dy)+(\ed_1 dx_1-\ed_3 dx)\wedge dy_n\\
		&+dx_n\wedge(\ed_1 dy_1-\ed_3dy)\\
		=\frac{\ed_3^2}{\ve_2}(dx_n&\wedge dy_2+dx_2\wedge dy_n)-\frac{\ed_3^2}{\ve_2}(dx_n\wedge dy_2+dx_2\wedge dy_n)=0,
	\end{split}
\end{equation*}
as needed.
\end{proof}

\begin{remark}
{\rm It was pointed out by A. Izosimov that the 2-form $\omega$ is also well defined on the space $\widetilde{\mathcal X}_{n,{\bf S},{\mathcal M}}$ of twisted $n$-gons with monodromy ${\mathcal M}$.
}
\end{remark}

\subsection{Additional integrals} \label{subsec:add}

The next theorem, providing two additional integrals of the $c$-relation, is a discrete version of Proposition 3.4 in \cite{Ta18}, concerning a continuous version of the $c$-relation on centroaffine curves. It is also an analog of Theorem 16 in \cite{AFIT}.   

\begin{theorem} \label{thm:extra}
1) The functions $I,J,K$ are integrals of the $c$-relations for all values of $c$ and of the polygon recutting.\\
2) The integral $4IK-J^2$ descends to the moduli space ${\mathcal X}_{n,{\bf S}}$. 
\end{theorem}

\begin{proof}
We use the notations from the proof of Theorem \ref{thm:pres}. 

One has
$$
x_i v_i - y_i u_i = c, \
x_i y_{i+1} - y_i x_{i+1} = s_{2i+1} = u_i v_{i+1} - v_i u_{i+1}
$$
for all $i$. Hence $v_i=(c+y_iu_i)/x_i$ and
$$
\frac{u_i (c+y_{i+1}u_{i+1})}{x_{i+1}} - \frac{u_{i+1}(c+y_iu_i)}{x_i} = s_{2i+1}.
$$
It follows that 
$$
c(x_iu_i-x_{i+1}u_{i+1}) + u_iu_{i+1} (x_iy_{i+1}-y_ix_{i+1}) = s_{2i+1} x_ix_{i+1}
$$
or
$$
c[x_iu_i-x_{i+1}u_{i+1}] + s_{2i+1}  u_iu_{i+1} = s_{2i+1} x_ix_{i+1}.
$$
Taking sum over $i$ eliminates the first summand on the left hand side and shows that $I$ is an integral.

Next we show that $I$ is also an integral for the polygon recutting. Indeed, 
\[I(\mathbf{P})-I(\rec_{1}(\mathbf{P}))=(x_1x_2\ed_3+x_nx_1 \ed_1)-(x x_2 \ed_1+x_n x \ed_3)=\]\[=x_2(x_1\ed_3-x\ed_1)+x_n(x_1\ed_1-x\ed_3)=0.\]
Therefore $I(\mathbf{P})=I(\rec_{1}(\mathbf{P}))=I(\rec_{2}\rec_1(\mathbf{P}))=...=\rec(\mathbf{P})$.    

One also has
\begin{equation} \label{eq:slact}
\begin{split}
e(I)=0, e(J)=2I, e(K)=J,\ &h(I)=2I, h(J)=0, h(K)=-2K,\\
 &f(I)=J, f(J)=2K, f(K)=0.
\end{split}
\end{equation}
Since the $c$-relations and the  recutting commute with $SL(2,\R)$, it follows that $J$ and $K$ are also integrals.

Finally, (\ref{eq:slact}) imply that $4IK-J^2$ is an $sl(2,\R)$-invariant function. This proves the last claim.
\end{proof}

The space spanned by $I,J,K$ is the irreducible 3-dimensional representation of the Lie algebra $sl(2,\R)$, the symmetric square of its standard 2-dimensional representation, or the coadjoint representation. The map ${\mathcal X}_{n,{\bf S}}\to \R^3$, whose components are the functions $I,J,K$, is the moment map of the Hamiltonian action of $sl(2,\R)$ on ${\mathcal X}_{n,{\bf S}}$.

\begin{remark}
{\rm One has
$$
4IK-J^2 =- \sum_{k,l} s_{2k+1} s_{2l+1}( [P_k,P_l] [P_{k+1} P_{l+1}]-[P_k,P_{l+1}] [P_{k+1},P_l]).
$$
}
\end{remark}

\subsection{Center of a polygon} \label{subsec:center}

Define the {\it center} of a polygon ${\bf P}$ as the quadratic form on $\R^2$ 
$$
C({\bf P})=I({\bf P}) x^2 - J({\bf P}) xy + K({\bf P}) y^2.
$$
The center is invariant under the $c$-relation and the polygon recutting, and it conjugates the diagonal action of $SL(2,\R)$ on polygons and on its action  on quadratic forms. In this section we present properties of the center, somewhat analogous to those of the circumcenter of mass, see \cite{TT14}. 

The center is additive in the following sense. 

\begin{lemma} \label{lm:addit}
Let ${\bf P}=(P_1,\ldots,P_k,\ldots,P_n)$. Cut ${\bf P}$ into two polygons ${\bf P}_1=(P_1,\ldots,P_k)$ and ${\bf P}_2=(P_1,P_k,\ldots,P_n)$. 
Then $C({\bf P})=C({\bf P}_1) + C({\bf P}_2)$.
\end{lemma} 

\begin{proof}
Each of the three components of the sum $C({\bf P}_1) + C({\bf P}_2)$ contains all the terms of the respective component of $C({\bf P})$ plus the additional terms (two or four), that appear due to the cut $P_1P_k$. Since the ``side area" $[P_1P_k]$ changes the sign when the orientation of the side is reversed, these additional terms cancel pairwise.
\end{proof}

According to the preceding lemma, the calculation of the center of a polygon reduces to that of a triangle. The next result gives a geometrical interpretation to the center of a triangle.

\begin{lemma} \label{lm:centtri}
A triangle ${\bf P}=(P_1P_2P_3)$ admits a unique circumscribed central conic given by the equation $ax^2-bxy+cy^2=1$. The  center of the triangle ${\bf P}$ is $2A(ax^2-bxy+cy^2)$, where $A$ is the oriented area of ${\bf P}$.
\end{lemma}

\begin{proof}
Let $P_i=(x_i,y_i),\ i=1,2,3$. We recall that the sides to not pass through the origin.

To find the circumscribed central conic one needs to solve the linear system $M (a,-b,c)^T=(1,1,1)^T$ where
$$
M= \begin{pmatrix}
	x_1^2&x_1y_1&y_1^2\\
	x_2^2&x_2y_2&y_2^2\\
	x_3^2&x_3y_3&y_3^2
\end{pmatrix}.
$$
One has $\det M = \ed_1\ed_3\ed_5$. Denote by $N$ the cofactor matrix of $M$. Then
$$
M^{-1} = \frac{1}{\ed_1\ed_3\ed_5} N,
$$
and hence 
$$
(a,-b,c)^T = \frac{1}{\ed_1\ed_3\ed_5} N (1,1,1)^T.
$$
Note that $\ed_1\ed_3\ed_5$ is twice the oriented area of the triangle.

On the other hand, one notices that $N (1,1,1)^T = (K,-J,I)^T$. This implies the result. 
\end{proof}

\begin{lemma} \label{lm:centbut}
The center of a centroaffine butterfly is the origin.
\end{lemma}

\begin{figure}[ht]
\centering
\includegraphics[width=.53\textwidth]{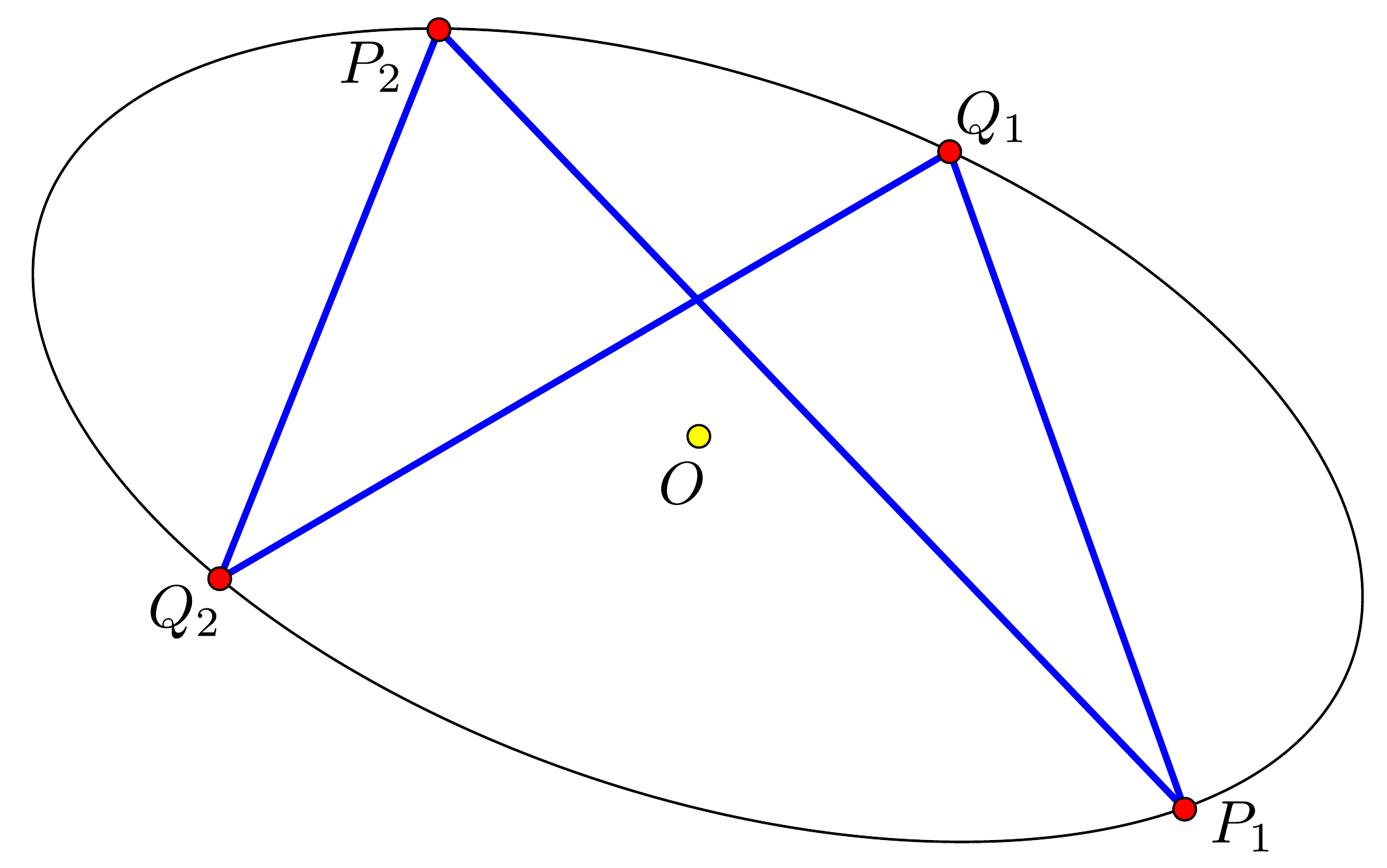}
\caption{Triangles $P_1P_2Q_2$ and $P_1Q_1Q_2$ share the circumcribed central conic.}
\label{Center}
\end{figure}

\begin{proof} An affine reflection in a line through the origin interchanges points $P_1$ with $Q_2$, and $P_2$ with $Q_1$, see Figure \ref{Center}. This reflection preserves the central conic circumscribed about triangle $P_1P_2Q_2$, therefore triangles $P_1P_2Q_2$ and $P_1Q_1Q_2$ share the circumcribed central conic. Now the result follows from Lemmas 
\ref{lm:addit} and \ref{lm:centtri}.
\end{proof}

Lemmas \ref{lm:addit} and \ref{lm:centbut} provide an alternative proof that the center is invariant under polygon recutting.

\section{Small-gons} \label{sec:small}

\subsection{Triangles}
In this section we investigate closed triangles. 
 
\begin{figure}[ht]
\centering
\includegraphics[width=.6\textwidth]{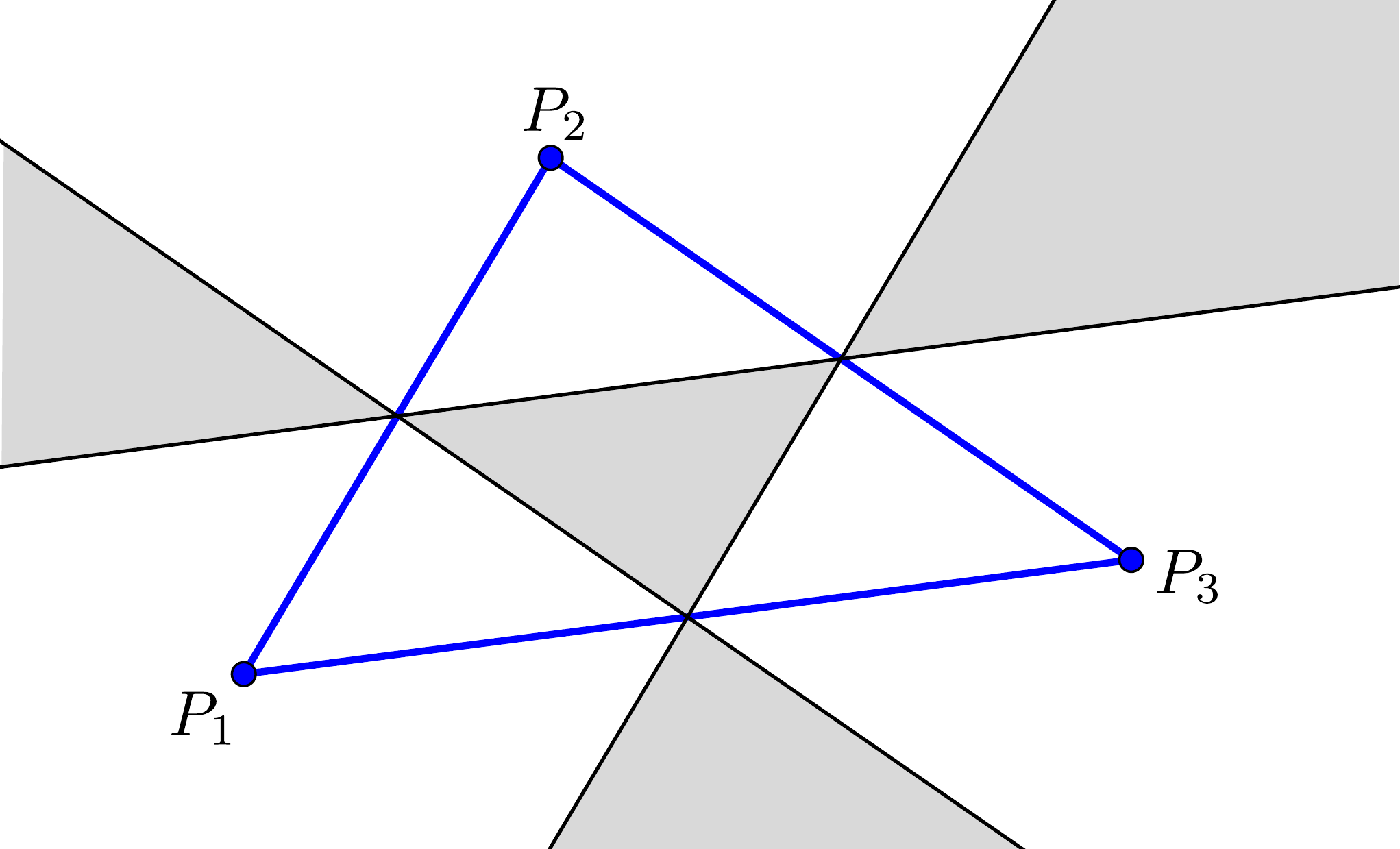}
\caption{The shaded regions are where the origin should be located for the linear map $M$ to be elliptic.}
\label{middle}
\end{figure}

\begin{theorem} \label{thm:tri}
1) A triangle admits a $c$-related triangle if and only if 
$$
c^2 (\ed_1+\ed_3+\ed_5)(\ed_1+\ed_3-\ed_5)(\ed_3+\ed_5-\ed_1)(\ed_5+\ed_1-\ed_3) \leq 4(\ed_1\ed_3\ed_5)^2.
$$
No triangles have infinitely many $c$-related ones for any $c\ne 0$.\\
2) Two $c$-related triangles are $SL(2,\R)$-equivalent. \\
3) The linear transformation $M$ that relates them and that defines the dynamics is elliptic if and only if
\begin{equation} \label{eq:type}
 (\ed_1+\ed_3+\ed_5)(\ed_1+\ed_3-\ed_5)(\ed_3+\ed_5-\ed_1)(\ed_5+\ed_1-\ed_3)>0.
\end{equation}
$M$ is parabolic if and only if the origin is located on the lines that bisects two sides of the triangle, see Figure \ref{middle}.\\
4) Let triangle $A'B'C'$ be the recutting of triangle $ABC$ done in the order $A \to B \to C$, and let $A''B''C''$ be the second iteration of this recutting. Then there exists a transformation $M\in SL(2,\R)$  that takes $ABC$ to $A''B''C''$.
\end{theorem}

\begin{proof} A triangle  is uniquely determined, modulo $SL(2,\R)$, by the areas $\ed_1,\ed_3,\ed_5$. These numbers are preserved by the $c$-relation, proving the second claim.

Let ${\bf Q}=M {\bf P}$ where $M=\begin{pmatrix}
	m&n\\k&l
\end{pmatrix}\in SL(2,R)$, and let $P_j=(x_j,y_j)$. Then the relation the ${\bf P}\crel {\bf Q}$ implies
$$\begin{cases}
	(m-l) x_1y_1-kx_1^2+ny_1^2=c\\
	(m-l) x_2y_2-kx_2^2+ny_2^2=c\\
	(m-l)x_3y_3-kx_3^2+ny_3^2=c.
\end{cases}$$
Here is the solution:
\[\begin{pmatrix}
	m-l\\
	k\\
	n
\end{pmatrix}=c\begin{pmatrix}
x_1y_1&-x_1^2&y_1^2\\
x_2y_2&-x_2^2&y_2^2\\
x_3y_3&-x_3^2&y_3^2
\end{pmatrix}^{-1}\begin{pmatrix}
1\\1\\1
\end{pmatrix}
.\]

For $M$ to exist, one needs the relation $ml-kn=1$ to hold. Since $m-l,k,n$ are already determined, this
reduces to a quadratic equation on $m$ that has real roots  if and only if 
$$
(m-l)^2+4kn +4 \ge 0.
$$

One has a remarkable identity:
\[(m-l)^2+4kn=-c^2\frac{ (\ed_1+\ed_3+\ed_5)(\ed_1+\ed_3-\ed_5)(\ed_3+\ed_5-\ed_1)(\ed_5+\ed_1-\ed_3)}{(\ed_1\ed_3\ed_5)^2},
\]
that we verified using Mathematica. This  implies the first claim of the theorem.
Furthermore,  $M$ is elliptic if and only if 
$$
\mathrm{tr}(M)^2-4\det (M)=(m+l)^2-4 < 0
$$ or, which is equivalent, $(m-l)^2+4kn < 0$. 
This implies the third claim.

 The right hand side vanishes when  $\ed_1=\ed_3+\ed_5$ or its cyclic permutation, that is,  
 the origin lies on one of the three middle lines of the triangle. These lines separate the elliptic and hyperbolic regions.
 
 For the last claim, one has $[A,B]=[B',C']$ and $[B,C]=[A',B']$, see Figure \ref{recut3}. Since the total area is preserved by recutting, one also has $[C,A]=[C',A']$. 

\begin{figure}[ht]
\centering
\includegraphics[width=.65\textwidth]{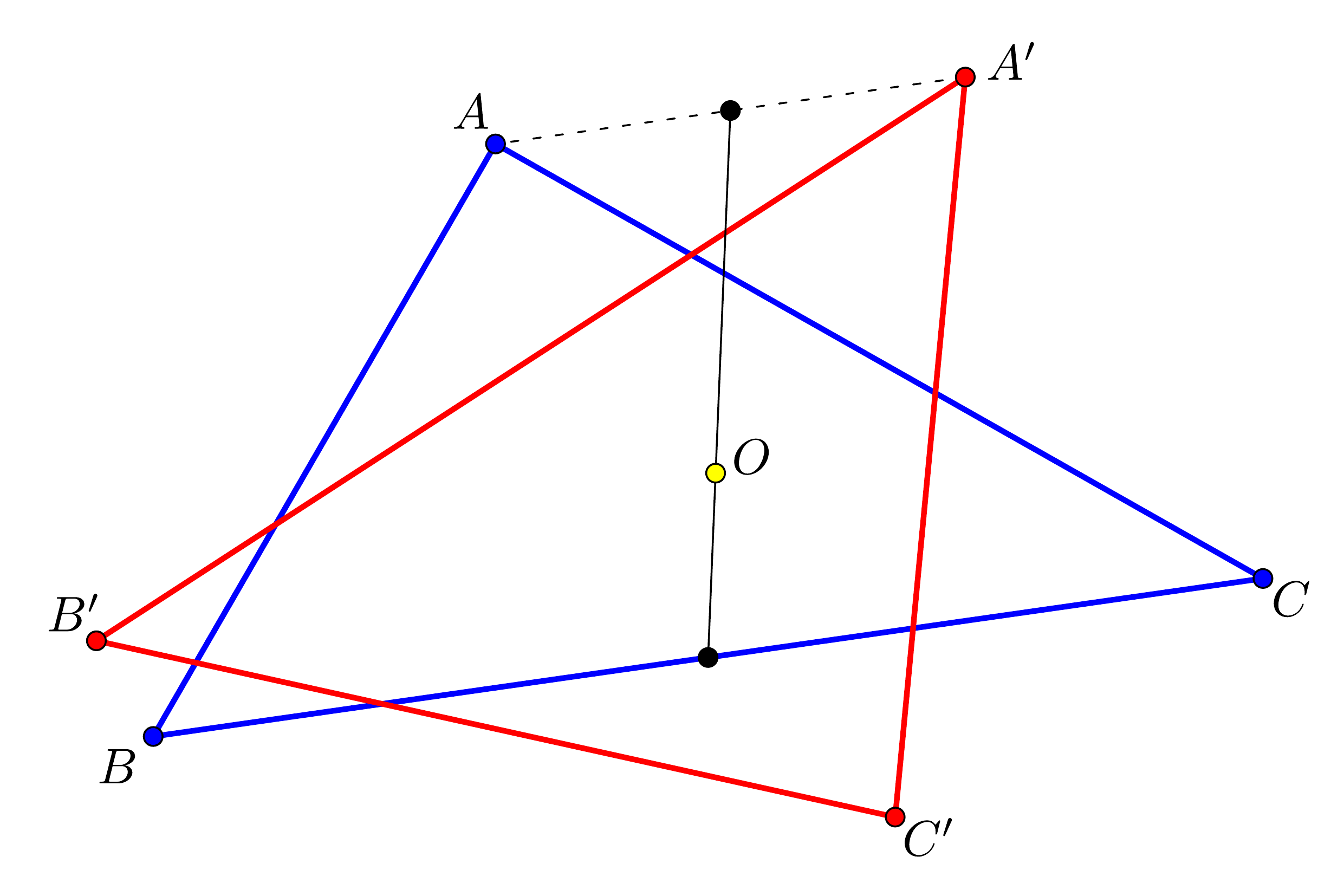}
\caption{Triangle recutting.}
\label{recut3}
\end{figure}

Repeating this argument going from  $A'B'C'$ to $A''B''C''$, we see that
$$
[A,B]=[A'',B''],\ [B,C]=[B'',C''],\ [C,A]=[C'',A''].
$$
Therefore the triangles $ABC$ and $A''B''C''$ are $SL(2,\R)$-equivalent.
 \end{proof}

 \begin{remark}
	{\rm The first claim if Theorem \ref{thm:tri} has the following geometric interpretation. Given a  triangle, there exists a central conic through its vertices. Assume that this conic is an ellipse and apply a transformation from $SL(2,\R)$ to make this ellipse into a circle of radius $R$. 
		Then a $c$-related triangle is also inscribed in this circle, and one has $c\le R^2$.
		
		Thus one expects the following identity to hold:
		\begin{equation} \label{eq:trig}
			R^4 (\ed_1+\ed_3+\ed_5)(\ed_1+\ed_3-\ed_5)(\ed_3+\ed_5-\ed_1)(\ed_5+\ed_1-\ed_3) = 4(\ed_1\ed_3\ed_5)^2.
		\end{equation}
		Let $\alpha, \beta, \gamma$ be the (signed) angles under which the sides of the triangle are seen from the origin. Then $\alpha + \beta +\gamma = 2\pi$,  
		$$
		\ed_1=R^2 \sin \alpha,\ \ed_3=R^2 \sin \beta, \ \ed_5=R^2 \sin \gamma,
		$$
		and (\ref{eq:trig}) becomes a true trigonometric identity
		\begin{equation*}
			\begin{split} 
				(\sin\alpha+\sin\beta+\sin\gamma)& (\sin\alpha+\sin\beta-\sin\gamma) (\sin\alpha-\sin\beta+\sin\gamma) \\
				&(-\sin\alpha+\sin\beta+\sin\gamma) = 4 \sin^2 \alpha \sin^2 \beta \sin^2 \gamma.
			\end{split}
		\end{equation*}
	}
\end{remark}

\subsection{Quadrilaterals} \label{subsec:quad}

Let us consider the dynamics of the $c$-relation on closed quadrilaterals. 
Let ${\mathbf P}$ be a quadrilateral, and assume that ${\mathcal L}_{{\bf P},c} \neq Id$.

\begin{figure}[ht]
\centering
\includegraphics[width=.47\textwidth]{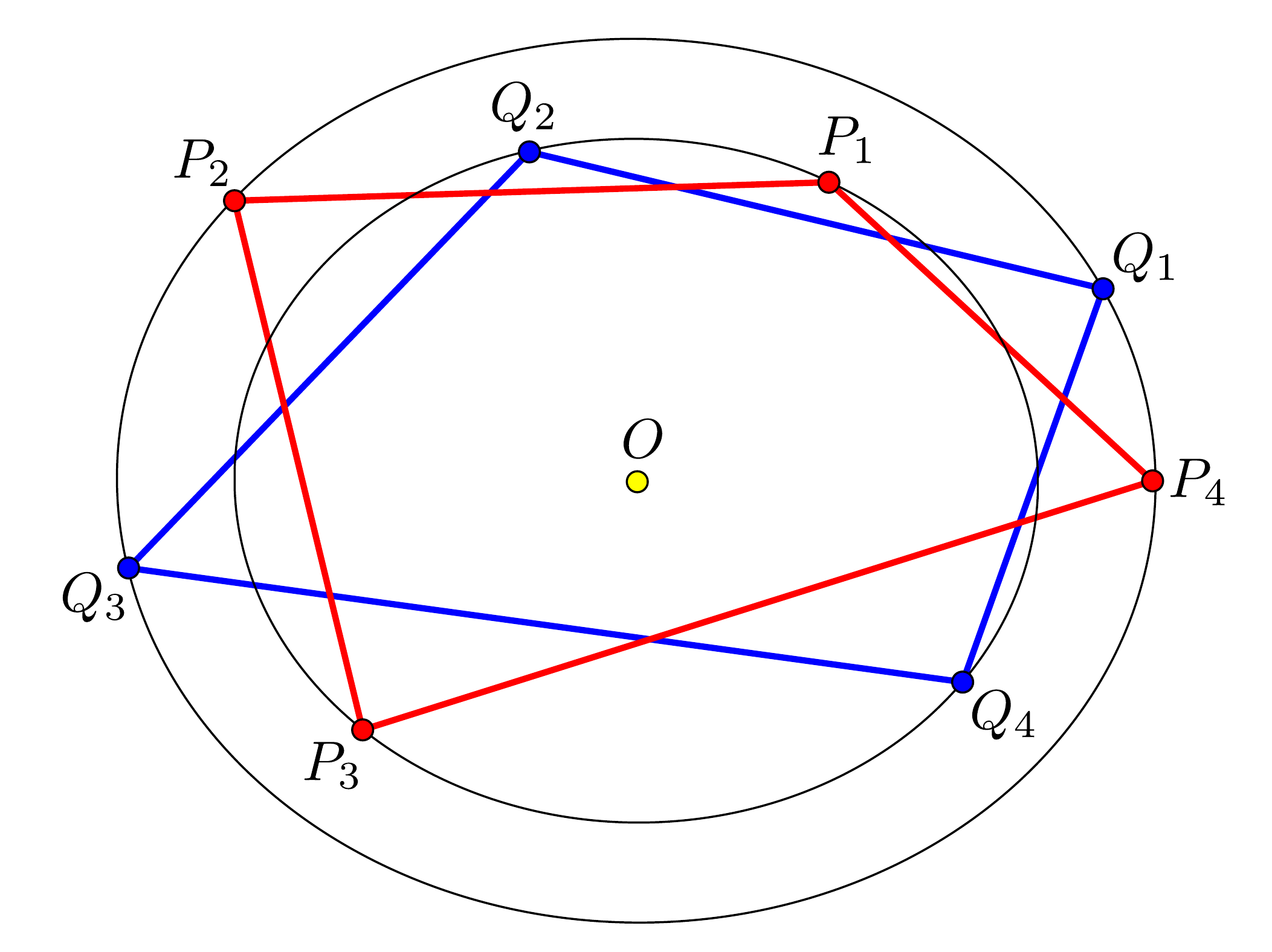}
\includegraphics[width=.47\textwidth]{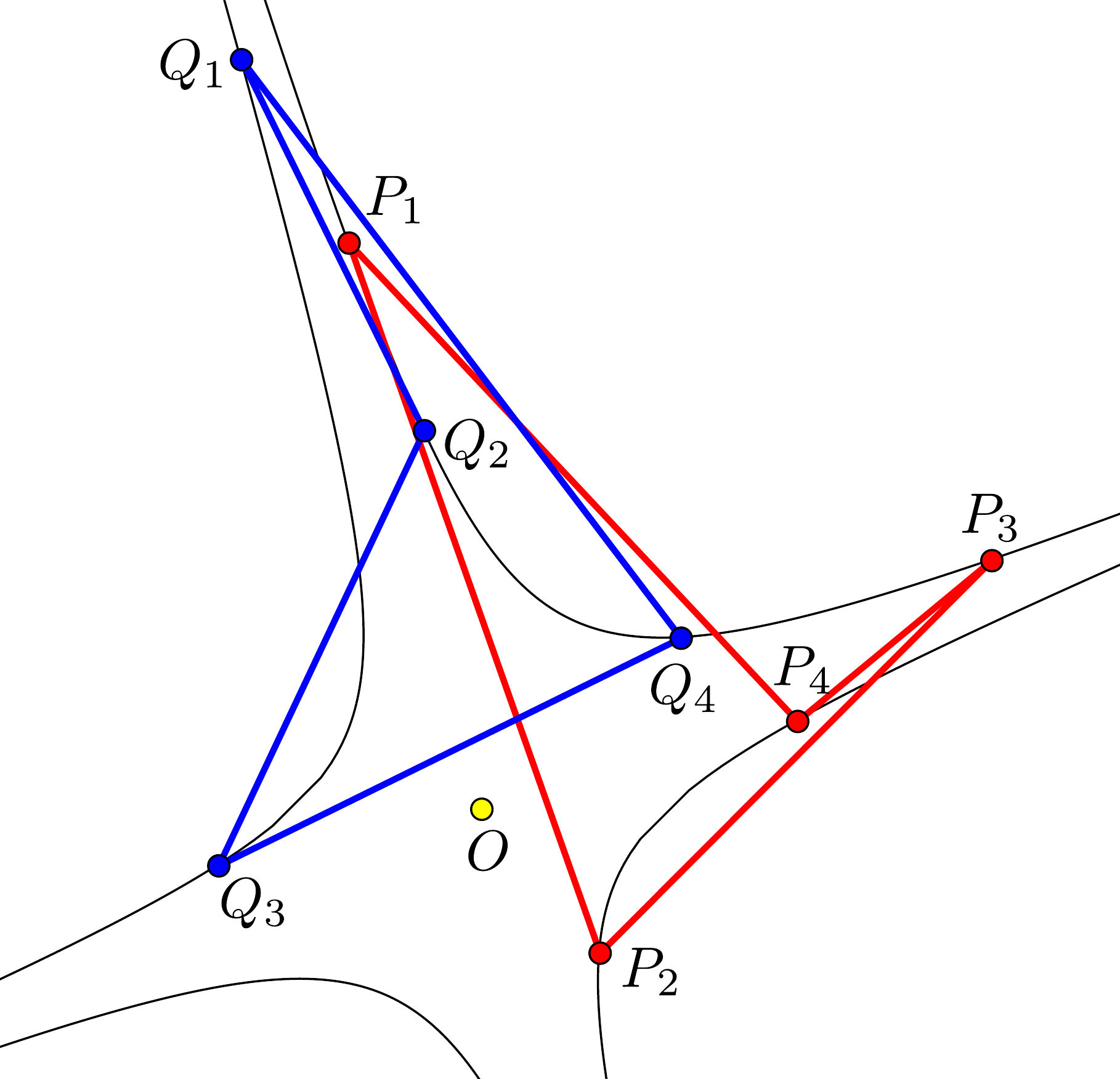}
\caption{Homothetic concentric conics containing the vertices of $c$-related quadrilaterals.}
\label{4conics}
\end{figure}

\begin{proposition} \label{prop:4dyn}
1) Let ${\bf P} \crel {\bf Q}$. Then there exist homothetic central conics ${\mathcal C_1}$ and ${\mathcal C_2}$ such that
$ P_1,Q_2,P_3,Q_4 \in {\mathcal C_1}$ and  $Q_1,P_2,Q_3,P_4 \in {\mathcal C_2}$, see Figure \ref{4conics}.\\
2) The conics in questions are ellipses if and only if 
$$
(\ed_1+ \ed_{3}+\ed_{5}+\ed_{7})(\ed_1+\ed_{3}-\ed_{5}-\ed_{7})(\ed_{3}+\ed_{5}-\ed_{7}-\ed_{1})(\ed_{5}+\ed_{7}-\ed_{1}-\ed_3)<0.
$$
\end{proposition}

\begin{proof}
As in Section \ref{subsec:constr}, one has four affine reflections $R_i$ whose composition is the identity map. Consider $R_i$ as an orientation reversing isometry of the hyperbolic plane, a reflection in a line $\ell_i$. Then $R_2 \circ R_1$ is an orientation preserving isometry. One has four cases: this isometry is elliptic, hyperbolic, parabolic, or the identity.

In the elliptic case, the isometry is a rotation about a point in $H^2$. Hence $R_4\circ R_3$ is a rotation about the same point, and therefore the four lines $\ell_i$ are concurrent at this point. It follows that  $R_i$ belong to  group $G\subset GL(2,\R)$ that is conjugated to $O(2)$, that is, the group generated by the rotations 
$$
\begin{pmatrix}
	\cos t&\sin t\\
	-\sin t&\cos t
\end{pmatrix}
\ \ {\rm and \ by}\ \ 
\begin{pmatrix}
	-1&0\\
	0&1
\end{pmatrix},
$$
and preserving the quadratic form $x^2+y^2$.
Thus $G$ preserves a positive-definite quadratic form whose circles are the homothetic ellipses preserved by the reflections $R_i$. Since 
$R_i$ swaps $P_i$ with $Q_{i+1}$ and $Q_i$ with $P_{i+1}$, we are done in this case.

In the hyperbolic case, the argument is similar. The isometries $R_i$ of $H^2$ are reflections in the lines $\ell_i$ that share a common perpendicular (and the lines are concurrent at a point of the projective plane outside of the absolute). In this case the argument is similar with 
the group $G$ being conjugated to $O(1,1)$ and generated by
$$
\begin{pmatrix}
	\cosh t&\sinh t\\
	\sinh t&\cosh t
\end{pmatrix},
\begin{pmatrix}
	-1&0\\
	0&1
\end{pmatrix},
\ \ {\rm and \ by}\ \ 
\begin{pmatrix}
	-1&0\\
	0&-1
\end{pmatrix},
$$
preserving the quadratic form $x^2-y^2$.
 Thus $G$ preserves a non-degenerate sign-indefinite quadratic form whose level curves are the desired homothetic hyperbolas (note that one of the level curves is singular: it is a pair of lines intersecting at the origin). 

\begin{figure}[ht]
\centering
\includegraphics[width=.5\textwidth]{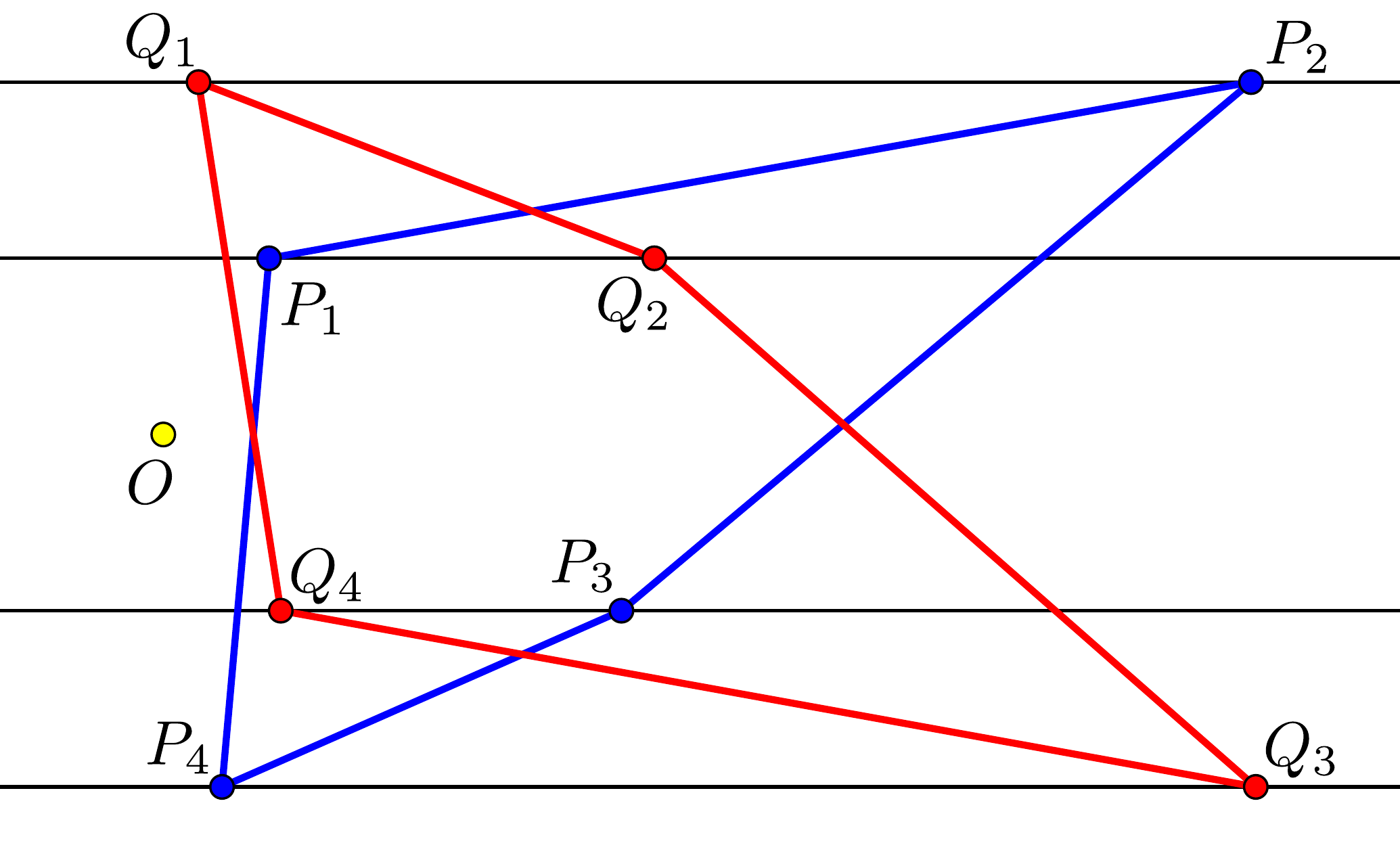}
\caption{The parabolic case.}
\label{lines}
\end{figure}

In the parabolic case, the conics degenerate to pairs of origin-symmetric parallel lines, see Figure \ref{lines}. Another degenerate case is when two opposite vertices of ${\mathbf P}$ are collinear, this happens in the hyperbolic case when the zero level curve of the sign-indefinite quadratic form is a pair of intersecting lines; see Figure \ref{lines} on the right.

Finally, in the case of the identity, one has $P_1 = (R_2 \circ R_1)(P_1)=P_3$, contradicting the non-degeneracy of the quadrilateral. 

We write the conics in the format $\langle P, MP\rangle$=const, where $P=(x,y)$ and $M=\begin{pmatrix}
		m&n\\n&k
	\end{pmatrix}$.
	The homogeneous equation for the matrix elements 
	\[\begin{cases}
		m(x_1^2-x_3^2)+2n(x_1y_1-x_3y_3)+k(y_1^2-y_3^2)=0\\
		m(x_2^2-x_4^2)+2n(x_2y_2-x_4y_4)+k(y_4^2-y_4^2)=0
	\end{cases}\]
has the solution 
\[(m:n:k)=\begin{bmatrix}
	y_1y_2 \ed_{3}+y_2y_3 \ed_{5}+y_3y_4 \ed_{7}+y_4y_1 \ed_{9}\\
-(1/2)((x_1^2-x_3^2)(y_2^2-y_4^2)-(x_2^2-x_4^2)(y_1^2-y_3^2))\\
x_1x_2 \ed_{3}+x_2x_3 \ed_{5}+x_3x_4 \ed_{7}+x_4x_1 \ed_{9}
\end{bmatrix}\]
A calculation shows that
\begin{equation*}
\begin{split}
&n^2-mk\\
&=(\ed_1+ \ed_{3}+\ed_{5}+\ed_{7})(\ed_1+\ed_{3}-\ed_{5}-\ed_{7})(\ed_{3}+\ed_{5}-\ed_{7}-\ed_{1})(\ed_{5}+\ed_{7}-\ed_{1}-\ed_3),
\end{split}
\end{equation*}
which implies the second result.	
\end{proof}

Arguing as in the preceding section, Proposition \ref{prop:4dyn} has the following corollary.

\begin{corollary} \label{cor:recut4}
	Let $\mathbf{Q}$ be the recutting of the quadrilateral $\mathbf{P}$. 
	Then the odd vertices of ${\bf Q}$ lie on the same central conic as the odd vertices of ${\bf P}$, and the even vertices of ${\bf Q}$ lie on the same homothetic central conic as the even vertices of ${\bf P}$.
\end{corollary}


Let ${\bf P}$ be a quadrilateral with coordinates $(s_1,s_3,s_5,s_7)$ and $(v_2,v_4)$. 

\begin{theorem} \label{thm:quad}
1) A quadrilateral admits a $c$-related quadrilateral if and only if
\begin{equation} \label{eq:cond4}
\begin{split}
c^2 (s_1 + s_3 - s_5 - s_7) &(s_1 - s_3 + s_5 - s_7) (s_1 - s_3 - s_5 + s_7) \\(s_1 + s_3 + s_5 + s_7)
  &\ge  4 (s_3 s_5 - s_1 s_7)  (s_3 s_7-s_1s_5) (s_1 s_3 - s_5 s_7).
\end{split}
\end{equation}
This condition is symmetric in $\{s_1,s_3,s_5,s_7\}$, and it has solution in $c$ for every $S$.\\
2) The second iteration of the $c$-transformation of ${\bf P}$ is $SL(2,\R)$-equivalent to ${\bf P}$.\\
3) The third iteration of the recutting of ${\bf P}$ is $SL(2,\R)$-equivalent to ${\bf P}$.
\end{theorem}

\begin{proof}
To prove the first statement, note that ${\bf P}$ is $SL(2,\R)$-equivalent to the following one:
$$
P_0=(1,0), P_1=(0,s_1), P_2=\left(-\frac{s_3}{s_1},v_2\right), P_3=\left(-\frac{v_4}{s_1},-s_7\right),
$$
where $v_2v_4=s_1s_5-s_3s_7$ (the Ptolemy-Pl\"ucker relation).

Let $Q_0=(b,c)$. A calculation using (\ref{eq:invol}) shows that $Q_4=(b',c)$, and we need $b'=b$. This is a quadratic equation $b^2+ub+v=0$ on $b$ whose coefficients are given by the formulas
\begin{equation*}
\begin{split}
&u=\frac{c}{v_2}\left[\frac{s_1^2-s_3^2-s_5^2-s_7^2}{s_5s_7-s_1s_3} + \frac{2s_3s_5s_7}{s_1(s_5s_7-s_1s_3)}\right],\\
&v=\frac{(c^2-s_1^2)(s_3s_5-s_1s_7)(s_3s_7-s_1s_5)}{v_2^2s_1^2(s_5s_7-s_1s_3)}.
\end{split}
\end{equation*}
One calculates the discriminant $D=u^2-4v$, and, after some cancellation, this results in (\ref{eq:cond4}) (we used Mathematica to clean-up the formulas).

One also has $D=c^2 P+4s_1^2 Q$, where 
\begin{equation*}
\begin{split}
&Q(s_1,s_3,s_5,s_7)=(s_5s_7-s_1s_3)(s_3s_7-s_1s_5)(s_3s_5-s_1s_7),\\
&P(s_1,s_3,s_5,s_7)=s_1^2(s_1^2-s_3^2-s_5^2-s_7^4)^2-Q(s_1,s_3,s_5,s_7).
\end{split}
\end{equation*}
Therefore one cannot have $P<0$ and $Q<0$ simultaneously, and this implies that (\ref{eq:cond4}) always has a solution.

The group of permutations of the four elements of the set $S$ is generated by the involutions $(13),(35),(57)$ that leave inequality (\ref{eq:cond4}) intact. 

The second statement of the theorem follows from Proposition \ref{prop:4dyn}: if ${\bf R} \crel {\bf Q}$, then the vertices of the quadrilateral ${\bf R}$ lie, alternating, on the same homothetic conics as those of ${\bf P}$, and the respective ``side areas" of these quadrilaterals are equal. This implies that ${\bf R}$ and ${\bf P}$ are $SL(2,\R)$-equivalent.

For the third statement, using the Ptolemy-Pl\"ucker relation, one calculates that  after the first recutting the coordinates become
$$
(s_1,s_5,s_7,s_3),  \left(v_2 \frac{s_1s_3-s_5s_7}{s_1s_5-s_3s_7},  v_4 \frac{s_1s_7-s_3s_5}{s_1s_3-s_5s_7}\right).
$$
Then after the second recutting we have
$$
(s_1,s_7,s_3,s_5), \left(v_2 \frac{s_1s_3-s_5s_7}{s_1s_7-s_3s_5},  v_4 \frac{s_1s_7-s_3s_5}{s_1s_5-s_3s_7}\right),
$$
so after the third recutting one obtains $(s_1,s_3,s_5,s_7), (v_2 ,  v_4)$, as claimed.
\end{proof} 

\subsection{Pentagons: invariant area form} \label{subsec:penta}

In this section we consider the moduli space ${\mathcal X}_{5,{\bf S}}$. This material is parallel to the one in Section 7.1.3 of \cite{AFIT}.

The space ${\mathcal X}_{5,{\bf S}}$ is 2-dimensional: the variables $\ve_{2i}$ satisfy the Ptolemy-Pl\"ucker relation  
$v_2v_4+v_8s_3=s_1s_5$ and its cyclic permutations.  We break the cyclic symmetry by setting $v_2=x, v_8=y$. Then
\begin{equation} \label{eq:Pl}
v_4=\frac{s_1s_5-s_3y}{x},\ v_6=\frac{s_5s_9-s_7x}{y},\ v_0=\frac{s_1s_7x+s_3s_9y-s_1s_5s_9}{xy}.
\end{equation}
Recall (Example \ref{ex:five})  that 
$$
K:=\frac{\ve_2}{\ed_1\ed_3} + \frac{\ve_4}{\ed_3\ed_5} + \frac{\ve_6}{\ed_5\ed_7} + \frac{\ve_8}{\ed_7\ed_9} + \frac{\ve_{10}}{\ed_9\ed_1},
$$
the only integral of the $c$-relation on the moduli space of closed pentagons. In terms of the $x,y$-coordinates, one has
$$
K=\frac{x}{s_1s_3}+\frac{y}{s_7s_9}+\frac{(s_1^2+s_3^2)}{s_1s_3x}+\frac{(s_7^2+s_9^2)}{s_7s_9y}-\frac{x}{s_5y}-\frac{y}{s_5x}-\frac{s_5}{xy}.
$$

Let 
$$
\omega=\frac{dx\wedge dy}{xy}.
$$
The origin of the area form on ${\mathcal X}_{5,{\bf S}}$ is in the theory of cluster algebras, and we do not dwell on it here.

\begin{theorem} \label{thm:pent}
The $c$-relation preserves the form $\omega$.
\end{theorem}

\begin{proof}
Recall the vector field $\xi$ from Section \ref{subsec:infin}. Using the formula from Theorem \ref{thm:field}, one has
\begin{equation} \label{eq:mixed}
\begin{split}
\dot x = x \left(\frac{v_4}{s_3s_5} - \frac{v_6}{s_5s_7} + \frac{v_8}{s_7s_9} - \frac{v_0}{s_9s_1}     \right) + \frac{s_3}{s_1}-\frac{s_1}{s_3},\\
\dot y = y \left(\frac{v_0}{s_9s_1}  - \frac{v_2}{s_1s_3} + \frac{v_4}{s_3s_5} - \frac{v_6}{s_5s_7}     \right) + \frac{s_3}{s_1}-\frac{s_1}{s_3}.
\end{split}
\end{equation}

We claim that $i_\xi \omega = dK$, that is, the integral $K$ is the Hamiltonian function of the vector field $\xi$. This claim is verified by a direct calculation, after substitution of the formulas (\ref{eq:Pl}) into (\ref{eq:mixed}).

Consider the $c$-relation as a transformation $T$. By the Bianchi permutability,  $\Phi_t$, the $t$-flow of the field $\xi$, commutes with $T$, and these maps preserve the level curves of the function $K$. Since $\Phi_t$ is symplectic, we have
$$
\Phi_t^* T^* (\omega) = T^* \Phi_t^* (\omega) = T^* (\omega),
$$ 
that is, $T^* (\omega)$ is also invariant under the flow of $\xi$. Hence $T^* (\omega)= H \omega$, where the function $H$ is an integral of $\Phi_t$ and of $T$ that has the same level curves as $K$. We want to show that $H\equiv 1$. 

Assume that a level curve $K=c$ is closed. Consider an infinitesimally close level curve $K=c+\eps$. Both curves are preserved by $T$, hence the area between them remains the same. On the other hand, this area is multiplied value of the function $H$ on the curve $K=c$. Hence this value equals 1, as needed.

\begin{figure}[ht]
	\centering
	\includegraphics[width=0.4\textwidth]{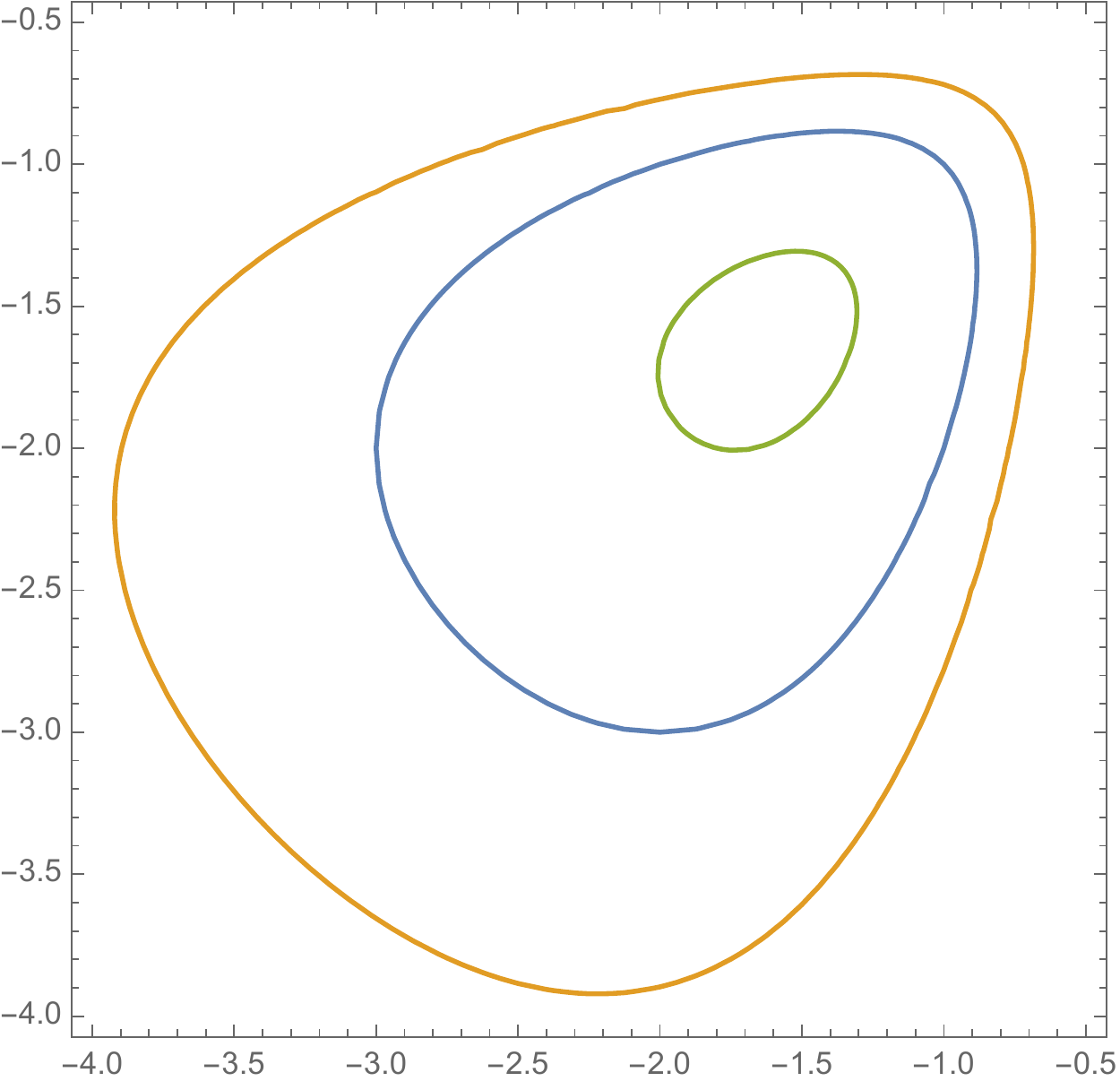}
	\caption{Three level curves $K=-10, -9, -8.2$ on ${\mathcal X}_{5,{\bf S}}$ with ${\bf S}=(1,1,1,1,1)$.}
	\label{curves}
\end{figure}

Thus $\omega$ is invariant under $T$ near local maxima or minima of $I$, see Figure \ref{curves}. The $c$-relation is an algebraic relation, and we can use analytic continuation to conclude that $c$-relation preserves $\omega$ everywhere.
\end{proof}

This theorem implies a Poncelet-style porism: if a level curve of the integral $K$ contains a periodic point of the $c$-relation, then every point of this curve is periodic with the same period.

\subsection{Pentagons: when the $c$-relation is not defined} \label{subsec:Fuchs}

Let ${\bf P}=(P_0,P_1,\ldots,P_{4})$ be a pentagon.   Using an appropriate $SL(2,{\mathbb R})$-transformation, we can make $P_0=(1,0)$ and $P_1=(0,s_1)$. Then
$$
P_2=\frac{(-s_3,s_1v_2)}{s_1},\ P_3=\frac{(s_3s_7-s_5v_4,s_1s_5s_9-s_1s_7v_4)}{v_2v_4-s_3s_9},\ P_4=\frac{(v_4,-s_1s_9)}{s_1}.
$$ 
To find a pentagon ${\bf Q}=(Q_0 Q_1,\ldots,Q_4)\ c$-related to $\bf P$, we first put $Q_0=(b,c)$ (with some unknown $b$) and compute the remaining vertices using the formulas from Section \ref{subsec:maps}. 

Our goal is to find a value of $b$ such that $Q_5=Q_0$; but for an arbitrary $b$ we will have $Q_5=(b',c)$ with some $b'$ which may be different from $b$ (because $[Q_5,P_5]=[Q_5,P_0]$ must be equal to $c$). This $b'$ will depend on $b,c$, and all $s,v$, and it is not hard to see that it will be, actually, a {\it quadratic} polynomial in $b$ with coefficients depending on $c,s,$ and $v$. Then the equation  $b'=b$ is quadratic with respect to $b$. 
(For $n=4$, a similar equation was explicitly calculated in the proof of Theorem \ref{thm:quad}.)

We were able to make these calculations, but the result looks depressive, and we do not present it here. Actually, we are more interested in the {\it discriminant} of this quadratic equation. This discriminant $D$ depends on $c,s_1,s_3,s_5,s_7,s_9$, and $v_2,v_4$, but in reality its dependence on $v_2,v_4$ may be reduced to the dependence on 
$$
K=\frac{v_2}{s_1s_3}+\frac{v_4}{s_3s_5}+\frac{v_6}{s_5s_7}+\frac{v_8}{s_7s_9}+\frac{v_{10}}{s_9s_1}.
$$
Moreover, $D$ turns out to be also a {\it quadratic} function of $K$ with coefficients depending on $c$, and $s_1,s_3,s_5,s_7,s_9$. An explicit expression for this function $D=D(K)$ is less awkward, it can be derived from Propositions \ref{prop:F1} and \ref{prop:F2} below. 

What we really need is the discriminant of $D(K)$, for which we will use a weird notation ${\mathcal D}(D)$. Indeed, {\it if, for some $c,s_1,s_3,s_5,s_7,s_9$, ${\mathcal D}(D)\leq0$, then $c$-related pentagons exist for all pentagons with these $s_1,s_3,s_5,s_7,s_9$; if ${\mathcal D}(D)>0$, then the equation $D(K)=0$ has two different real roots $K_1$ and $K_2$ and no $c$-related pentagons exist for pentagons with $K$ between $K_1$ and $K_2$.}

Notice that both $D$ and ${\mathcal D}(D)$ are defined up to positive factors, which we ignore in the formulas below.

\begin{proposition} \label{prop:F1}
One has
$${\mathcal D}(D)=\prod_{j=1}^5(c^2-s_{2j-1}^2).$$
\end{proposition}

\begin{figure}[ht]
	\centering
	\includegraphics[width=0.65\textwidth]{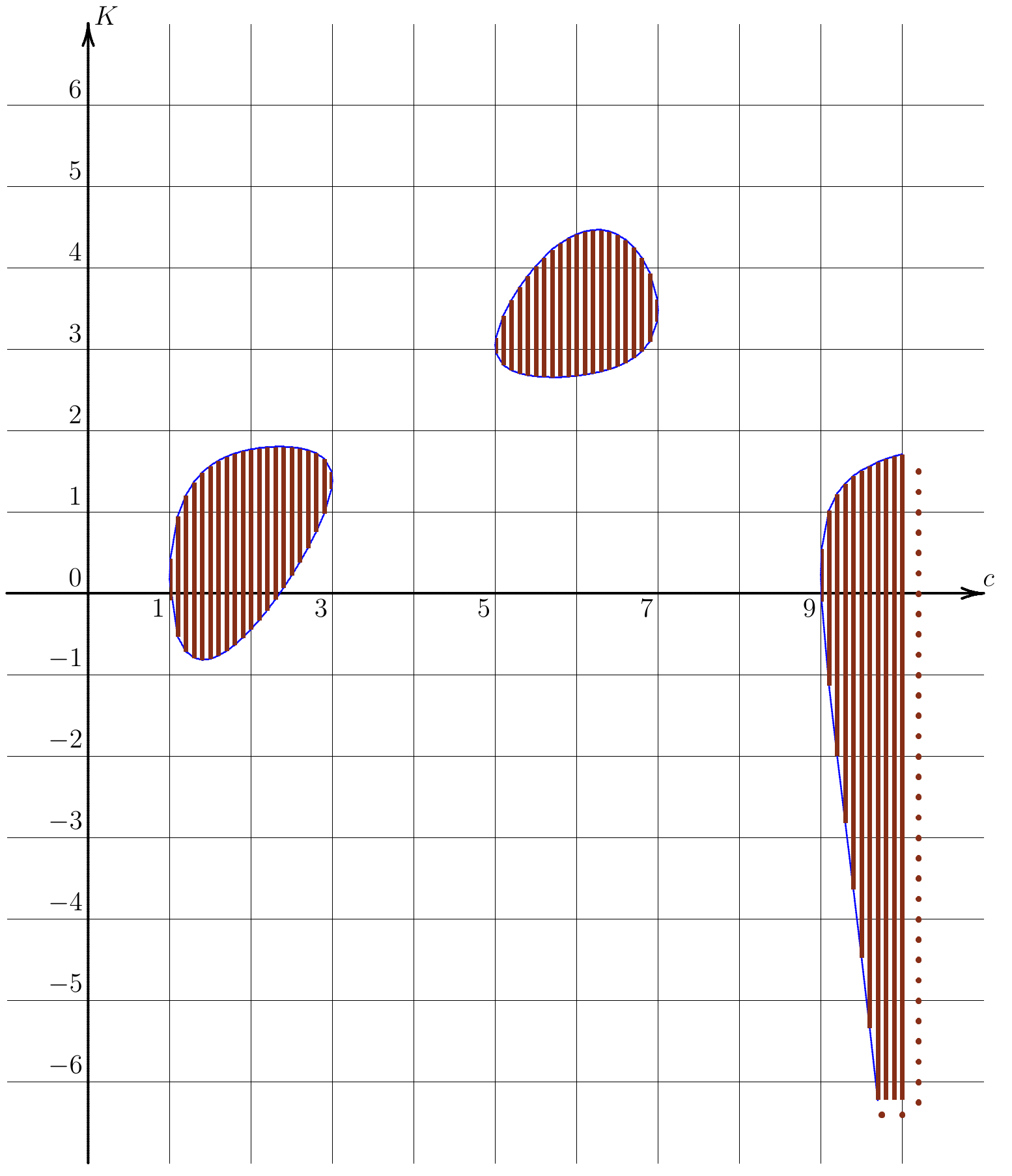}
	\caption{The zones where the $c$-relation is not defined; here ${\bf S}=\{1,3,5,7,9\}$. The horizontal and vertical axes are $c$ and $K$.}
	\label{discr}
\end{figure}

Thus, if $0<s_1<\ldots<s_9$ (this condition is not really restrictive, since all our results are not sensitive to permutations and sign changes of $s_1,\dots s_9$) then the relation may be undefined only if $s_1<|c|<s_3$, or $s_5<|c|<s_7$, or $|c|>s_9$. See Figure \ref{discr}.

\begin{proposition} \label{prop:F2}
The solutions of the equation $D(K)=0$ are
$$K=\frac{((\sum_js_{2j-1}^2)-2c^2)c^2}{\prod_js_{2j-1}}\pm\frac{2\sqrt{{\mathcal D}(D)}}{c\prod_js_{2j-1}}.$$
\end{proposition}

The proofs consist in tedious but explicit calculations. 

By the way, to check a reliable formula, after it has been obtained, we need, as a rule, to prove the equality between two polynomials of the same degree, and for this it is sufficient to verify the equality for a certain number of integral  of variables, which is an easy task for a computer program.

\end{document}